\documentclass[11pt]{amsart}

\usepackage{amsmath,amssymb,amsthm,amsrefs,a4wide}
\usepackage{latexsym}
\usepackage{indentfirst}
\usepackage{graphicx}
\usepackage{placeins}
\usepackage{booktabs}
\usepackage{algorithm}
\usepackage{algorithmic}
\usepackage{multirow}

\theoremstyle{plain}

\theoremstyle{definition}

\theoremstyle{remark}

\DeclareMathOperator{\argmin}{arg min}

\DeclareMathOperator{\cv}{cv}
\DeclareMathOperator{\rv}{rv}

\newcommand{\bbm}{\begin{bmatrix}}
\newcommand{\ebm}{\end{bmatrix}}
\newcommand{\R}{\mathbb{R}}
\newcommand{\C}{\mathbb{C}}
\newcommand{\D}{\mathbb{D}}

\newcommand{\Z}{\mathbb{Z}}

\newcommand{\p}{\partial}

\newcommand{\cst}{\text{const}}

\newcommand{\z}{\xi}
\renewcommand{\t}{\tau}

\newcommand{\hg}{\hat{g}}
\newcommand{\hG}{\hat{G}}

\newcommand{\dmax}{{d_{\text{max}}}}

\newcommand{\dt}{\mathrm{d}t}
\newcommand{\dtot}{\frac{\dt}{t}}
\newcommand{\dz}{\mathrm{d}z}
\newcommand{\dth}{\mathrm{d}\theta}

\begin{document}

\title[Pole recovery from noisy data on imaginary axis]{Pole recovery from noisy data on imaginary
  axis}

\author[]{Lexing Ying}

\address[Lexing Ying]{Department of Mathematics, Stanford University, Stanford, CA 94305}

\email{lexing@stanford.edu}

\thanks{The author thanks Lin Lin and Anil Damle for discussions on this topic.}


\begin{abstract}
  This note proposes an algorithm for identifying the poles and residues of a meromorphic function
  from its noisy values on the imaginary axis. The algorithm uses M\"{o}bius transform and Prony's
  method in the frequency domain. Numerical results are provided to demonstrate the performance of
  the algorithm.
\end{abstract}

\maketitle

\section{Introduction}\label{sec:intro}

\begin{figure}[h!]
  \centering
  \includegraphics[scale=0.3]{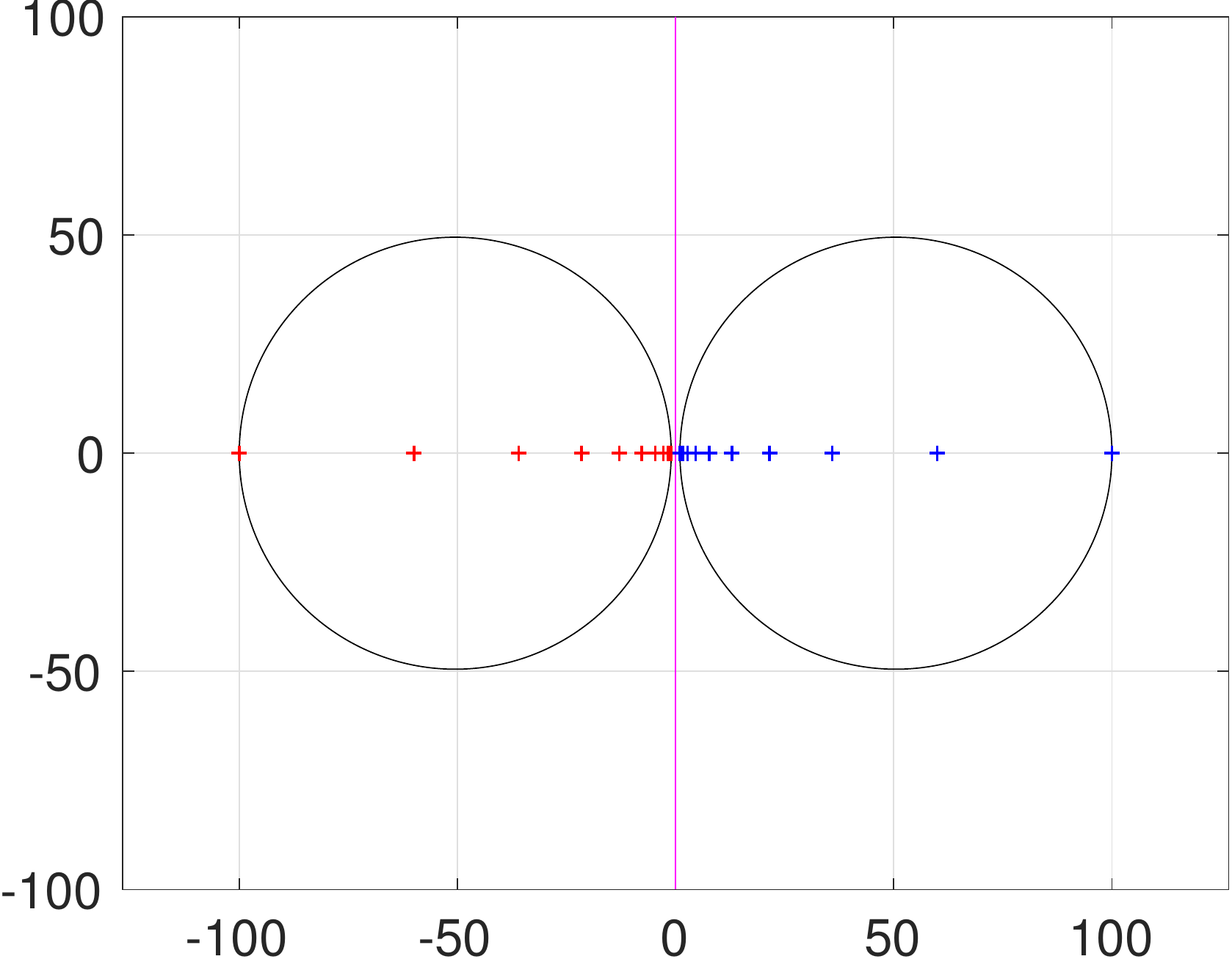}
  \caption{The unknown poles are inside the two circles. The algorithm can access the noisy function
    values along the imaginary axis.}
  \label{fig:problem}
\end{figure}

Let $g(z)$ be a meromorphic function of the form
\begin{equation}
  g(z) = \sum_{j=1}^{N_p} \frac{r_j}{\z_j-z}
  \label{eq:sclver}
\end{equation}
where the number of poles $N_p$, the pole locations $\{\z_j\}$, and residues $\{r_j\}$ are all
unknown, except that $\z_j$ are away from the imaginary axis $i \R$. The problem is to recover
$N_p$, $\{\z_j\}$ and $\{r_j\}$, given the noisy access of $g(z)$ along the imaginary axis $i
\R$. Two access models are particularly relevant: (1) the {\em random access model} where one can
get noisy values of $g(z)$ anywhere on $i\R$ and (2) the {\em Matsubara model} where one can only
get the noisy values of $g(z)$ at the Matusbara grid
\[
z_n =
\begin{cases}
  2n \frac{\pi i}{\beta},      & \text{for bosons},\\
  (2n+1) \frac{\pi i}{\beta},  & \text{for fermions}.
\end{cases}
\]
To make the problem numerically well-defined, we assume
\begin{itemize}
\item There exists constants $0<a<b$ such that the poles $\{\z_j\}$ reside within the union of the
  two disks of radius $\frac{b-a}{2}$ centered at $-\frac{b+a}{2}$ and $\frac{b+a}{2}$,
  respectively. See Figure \ref{fig:problem} for an illustration.
\end{itemize}
This assumption is quite natural because otherwise any algorithm is forced to sample extensively
along the imaginary axis towards infinity.

There is also a matrix-valued version of this problem, where
\begin{equation}
  G(z) = \sum_{j=1}^{N_p} \frac{R_j}{\z_j-z}
  \label{eq:matver}
\end{equation}
where $G(z)$ and $R_j$ are matrices of size $N_b\times N_b$. The task is then to recover $N_p$,
$\{\z_j\}$ and $\{R_j\}$. A particularly important special case is where $R_j = v_j v_j^*$ for some
$v_j \in \C^{N_b}$ \cite{mejuto2020efficient}.

This problem has many applications in scientific and engineering disciplines. One of the key
examples is the reconstruction of spectral density from Matsubara Green's function
\cite{bruus2004many}. This problem is highly related to a couple of other well-studied problems,
including rational function approximation and interpolation
\cite{berljafa2017rkfit,berrut2004barycentric,beylkin2009nonlinear,gustavsen1999rational,nakatsukasa2018aaa,wilber2021data},
Pade approximation \cite{gonnet2013robust}, contractive analytic continuation
\cite{fei2021nevanlinna,fei2021analytical}, approximation with exponential sums
\cite{beylkin2005approximation,potts2013parameter}, and hybridization fitting
\cite{mejuto2020efficient}.

Since this problem is quite ill-posed, a solution should be relatively robust to a reasonable level
of noise. The main content of this note is a simple algorithm based on conformal mapping and Prony's
method that naturally draws ideas from the references list above.


\section{Algorithm}\label{sec:algo}

\subsection{Continuous version}

Let us consider the scalar case \eqref{eq:sclver}. Below we describe the algorithm as if one can
manipulate continuous objects. The overall plan is to
\begin{itemize}
\item locate the poles in the left and right half plane separately using M\"{o}bius transform and
  Prony's method,
\item compute the residues using least square.
\end{itemize}

\begin{figure}[h!]
  \centering
  \includegraphics[scale=0.3]{tst0_A.pdf}\hspace{0.5in}
  \includegraphics[scale=0.3]{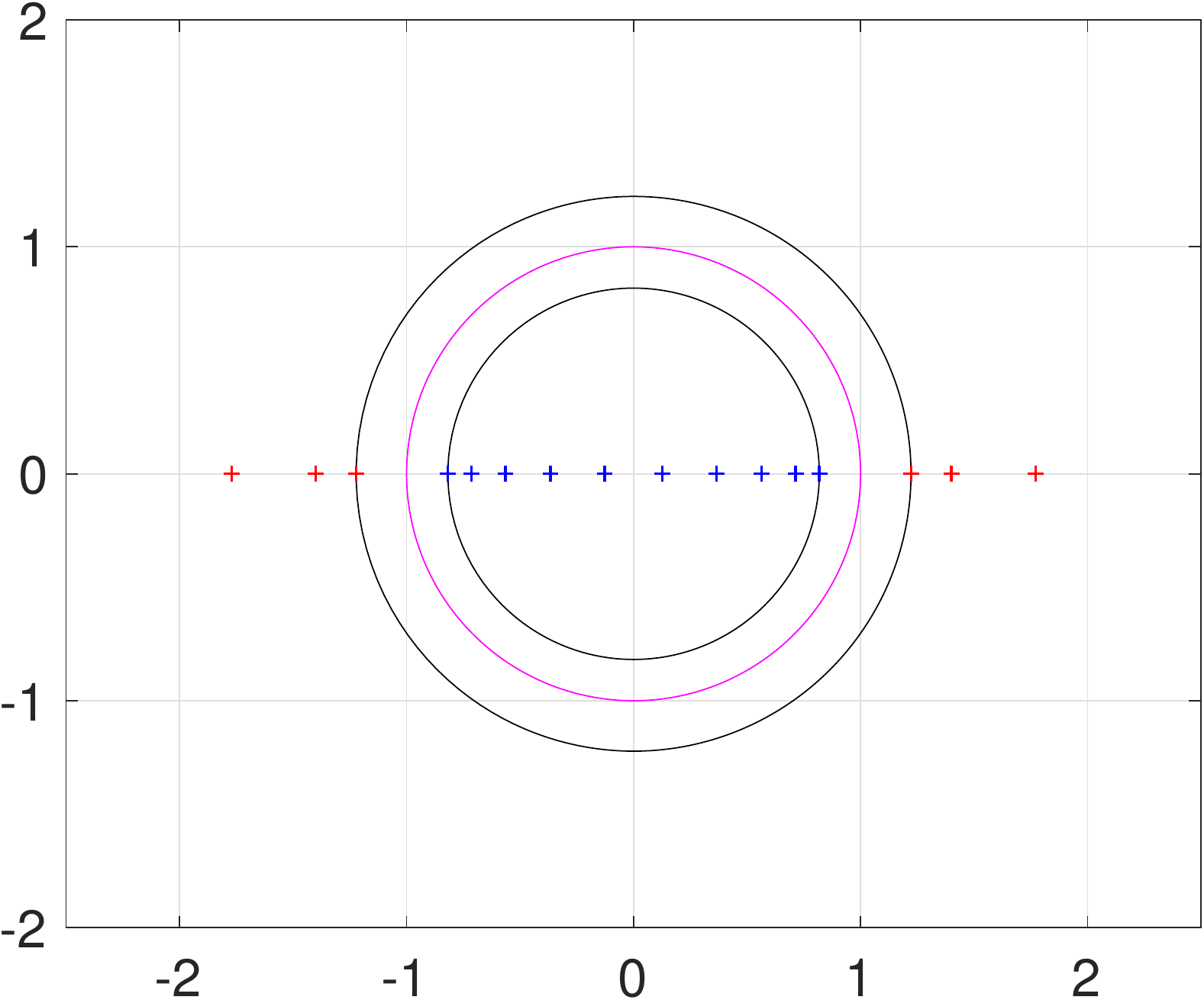}
  \caption{M\"{o}bius transform. Left: the $z$ plane. Right: the $t$ plane.}
  \label{fig:z2t}
\end{figure}

We introduce the following M\"{o}bius transform from $z\in\C$ to $t\in\C$
\begin{equation}
  t = \frac{z-\sqrt{ab}}{z+\sqrt{ab}}, \quad z = -\sqrt{ab} \frac{t+1}{t-1}.
  \label{eq:mobius}
\end{equation}
This transform maps
\begin{itemize}
\item the right half-plane $\C^+$ in $z$ to the interior of the unit disk $\D$ in $t$,
\item the left half-plane $\C^-$ in $z$ to the exterior of $\D$ in $t$,
\item the imaginary axis $i\R$ in $z$ to the unit circle in $t$, 
\item the two circles centered at $-\frac{b+a}{2}$ and $\frac{b+a}{2}$ in $z$ to two concentric
  circles with radius $\frac{\sqrt{b}-\sqrt{a}}{\sqrt{b}+\sqrt{a}}$ and
  $\frac{\sqrt{b}+\sqrt{a}}{\sqrt{b}-\sqrt{a}}$ in $t$ (see Figure \ref{fig:z2t} for an
  illustration).
\end{itemize}

The function $g(t)\equiv g(z(t))$ in the $t$ space also enjoys a pole representation
\[
g(z) = \sum_{j=1}^{N_p} \frac{w_j}{\t_j-t} + \cst
\]
for locations $\{\t_j\}$ and residues $\{w_j\}$. Since $\{\t_j\}$ are the images of the poles
$\{z_j\}$ under the M\"{o}bius transform, it is equivalent to locating $\{\t_j\}$.

Let us consider the integrals
\begin{equation}
  \frac{1}{2\pi i}\int_{\p\D} \frac{g(t)}{t^k} \dtot
  \label{eq:integral}
\end{equation}
for integer values of $k$. The integrals for negative and positive values of $k$ give information
about the poles inside $\D$ and the ones outside $\D$, respectively. For any $k\le -1$,
\[
\begin{aligned}
  & \frac{1}{2\pi i}\int_{\p\D}\frac{g(t)}{t^k}\dtot  = \frac{1}{2\pi i}\int_{\p\D}\left( \sum_{|\t_j|<1} + \sum_{|\t_j|>1}\right) \frac{w_j}{\t_j-t} t^{|k|+1} \dt \\
  &= \frac{1}{2\pi i} \sum_{|\t_j|<1} w_j \int \frac{1}{\t_j-t} t^{|k|-1} \dt
  = \frac{1}{2\pi i} \sum_{|\t_j|<1} w_j t_j^{|k|-1} \int \frac{1}{\t_j-t} \dt = - \sum_{|\t_j|<1} w_j \t_j^{-(k+1)},
\end{aligned}
\]
where the second equality uses the fact $\frac{w_j}{\t_j-t}$ is analytic in $\D$ for $|\t_j|>1$ and
the third equality uses the fact that $\frac{\t_j^{|k|-1}-t^{|k|-1}}{\t_j-t}$ is a polynomial hence
analytic in $\D$. Hence the integrals \eqref{eq:integral} for $k\le -1$ contain information about
the poles inside $\D$.

For any $k\ge 1$,
\[
\begin{aligned}
  & \frac{1}{2\pi i}\int_{\p\D} \frac{g(t)}{t^k} \dtot  = \frac{1}{2\pi i}\int_{\p\D} \left( \sum_{|\t_j|<1} + \sum_{|\t_j|>1}\right) \frac{w_j}{\t_j-t} \cdot \frac{1}{t^{k+1}} \dt \\
  &= \frac{1}{2\pi i} \sum_{|\t_j|>1} w_j \int_{\p\D} \frac{1}{\t_j-t} \cdot \frac{1}{t^{k+1}} \dt =
  \frac{1}{2\pi i} \sum_{|\t_j|>1} w_j \int_{\p\D} \frac{1}{\t_j} \left(1+\frac{t}{t_j}+\cdots  \right) \frac{1}{t^{k+1}} \dt \\
  &=  \frac{1}{2\pi i} \sum_{|\t_j|>1} w_j \frac{1}{\t_j^{k+1}} \int_{\p\D} \frac{1}{t} \dt
  = \sum_{|\t_j|>1} w_j \t_j^{-(k+1)},
\end{aligned}
\]
where the second equality uses the fact that for $|\t_j|<1$ the product
$\frac{w_j}{\t_j-t}\cdot\frac{1}{t^{k+1}}$ is analytic outside $\D$ with at least quadratic decay,
the fourth step uses the fact that only the term with $t^{k}$ in the power expansion gives
non-zero contribution.  Hence the integrals \eqref{eq:integral} for $k\ge 1$
contain information about the poles outside $\D$.

Since the integral \eqref{eq:integral} is over the unit circle, it is closely related to the
Fourier transform of the function $g(\theta) \equiv g(e^{i\theta})$:
\begin{equation}
  \frac{1}{2\pi i}\int_{\p\D} \frac{g(t)}{t^k} \dtot = \frac{1}{2\pi i}\int_0^{2\pi} g(\theta) e^{-ik\theta}  i \dth
  = \frac{1}{2\pi} \int_0^{2\pi} g(\theta) e^{-ik\theta} \dth = \hg_k.
  \label{eq:Fourier}
\end{equation}

To recover the poles inside $\D$, we use Prony's method to the Fourier coefficients. From the
integrals with $k\le -1$, define the semi-infinite vector
\[
\hg_-
\equiv
\begin{bmatrix}
  \hg_{-1}\\
  \hg_{-2}\\
  \vdots
\end{bmatrix}
\equiv
\frac{1}{2\pi i}\int_{\p\D} g(t) 
\begin{bmatrix}
  t^0\\
  t^1\\
  \vdots
\end{bmatrix}
\dt
\equiv
\begin{bmatrix}
  -\sum_{|\t_j|<1} w_j \t_j^0\\
  -\sum_{|\t_j|<1} w_j \t_j^1\\
  \vdots
\end{bmatrix}
\]
Let us define $S$ to be the shift operator that shifts the semi-infinite vector upward (and drops the
first element). For any $\t_j$ with $|\t_j|<1$, 
\[
S
\begin{bmatrix}
  \t_j^0\\
  \t_j^1\\
  \vdots
\end{bmatrix}
=
\begin{bmatrix}
  \t_j^1\\
  \t_j^2\\
  \vdots
\end{bmatrix},
\quad\text{i.e.,}\quad
(S-\t_j)
\begin{bmatrix}
  \t_j^0\\
  \t_j^1\\
  \vdots
\end{bmatrix}
= 0.
\]
Since the operators $S-\t_j$ all commute, 
\begin{equation}
\prod_{|\t_i|<1} \left(S- \t_i\right)
\begin{bmatrix}
  \t_j^0\\
  \t_j^1\\
  \vdots
\end{bmatrix}
= 0.
\label{eq:prodR}
\end{equation}
Since $\hg_-$ is a linear combination of such semi-infinite vectors, 
\[
\prod_{|\t_i|<1} \left(S- \t_i\right) \hg_- = 0.
\]
Suppose that $\prod_{|\t_i|<1} \left(t-\t_i\right) = p_0 t^0 + \cdots + p_d t^d$ with coefficients
$p_i$, where the degree $d$ is equal to the number of poles in $\D$. Then \eqref{eq:prodR} becomes
\begin{equation}
p_0 (S^0 \hg_-) + \cdots + p_d (S^d \hg_-) = 0,
\quad\text{i.e.,}\quad
\begin{bmatrix}
  \hg_{-1} & \hg_{-2} & \cdots & \hg_{-(d+1)} \\
  \hg_{-2} & \hg_{-3} & \cdots & \hg_{-(d+2)} \\
  \vdots & \vdots & \vdots & \vdots
\end{bmatrix}
\begin{bmatrix}
  p_0\\
  \ldots\\
  p_d
\end{bmatrix}
= 0.
\label{eq:lsR}
\end{equation}
This implies that the number of poles in $\D$ is equal to the smallest value $d$ such that the
matrix in \eqref{eq:lsR} is rank deficient. In addition, $(p_0,\ldots,p_d)$ can be computed as a
non-zero vector in the null-space of this matrix. Once $(p_0,\ldots,p_d)$ are available, the roots
of
\[
p(t) = p_0 t^0 + \ldots p_d t^d
\]
are the poles $\{\t_j\}$ inside $\D$.

To recover the poles outside $\D$, we use again Prony's method but to a different part of the
Fourier coefficients. From the integrals with $k\ge 1$, define the semi-infinite vector
\[
\hg_+  \equiv
\begin{bmatrix}
  \hg_{1}\\
  \hg_{2}\\
  \vdots
\end{bmatrix}
\equiv
\frac{1}{2\pi i}\int_{\p\D} g(t) 
\begin{bmatrix}
  t^{-2}\\
  t^{-3}\\
  \vdots
\end{bmatrix}
\dt
\equiv
\begin{bmatrix}
  \sum_{|\t_j|>1} w_j \t_j^{-2}\\
  \sum_{|\t_j|>1} w_j \t_j^{-3}\\
  \vdots
\end{bmatrix}
\]
With the same shift operator $S$, for any $\t_j$ with $|\t_j|>1$
\[
S
\begin{bmatrix}
  \t_j^{-2}\\
  \t_j^{-3}\\
  \vdots
\end{bmatrix}
=
\begin{bmatrix}
  \t_j^{-3}\\
  \t_j^{-4}\\
  \vdots
\end{bmatrix},
\quad\text{i.e.,}\quad
(S-\t_j^{-1})
\begin{bmatrix}
  \t_j^{-2}\\
  \t_j^{-3}\\
  \vdots
\end{bmatrix}
= 0.
\]
Since the operators $S-\t_j^{-1}$ all commute, 
\begin{equation}
  \prod_{|\t_i|>1} (S- \t_i^{-1})
  \begin{bmatrix}
    \t_j^{-2}\\
    \t_j^{-3}\\
    \vdots
  \end{bmatrix}
  = 0.
\label{eq:prodL}
\end{equation}
Since $\hg_+$ is a linear combination of such semi-infinite vectors, 
\[
\prod_{|\t_i|>1} (S- \t_i^{-1}) \hg_+ = 0.
\]
Suppose that $\prod_{|\t_i|>1} (t-\t_i^{-1}) = p_0 t^0 + \cdots + p_d t^d$ with coefficients $p_i$,
where the degree $d$ is equal to the number of poles outside $\D$. Then \eqref{eq:prodL} becomes
\begin{equation}
  p_0 (S^0 \hg_+) + \cdots + p_d (S^d \hg_+) = 0,
  \quad\text{i.e.,}\quad
\begin{bmatrix}
  \hg_{1} & \hg_{2} & \cdots & \hg_{d+1}  \\
  \hg_{2} & \hg_{3} & \cdots & \hg_{d+2} \\
  \vdots & \vdots & \vdots & \vdots
\end{bmatrix}
\begin{bmatrix}
  p_0\\
  \ldots\\
  p_d
\end{bmatrix}
= 0.
\label{eq:lsL}
\end{equation}
This implies that the number of poles outside $\D$ is equal to the smallest value $d$ such that the
matrix in \eqref{eq:lsL} is rank deficient. As before, $(p_0,\ldots,p_d)$ can be computed as a
non-zero vector in the null-space of this matrix and the roots of
\[
p(t) = p_0 t^0 + \ldots p_d t^d
\]
are $\{\t_j^{-1}\}$. Taking inverse of these roots gives the poles $\{\t_j\}$ outside $\D$.

Once the poles inside and outside $\D$ in the $t$ plane are ready, we take the union and apply
\eqref{eq:mobius} to get the poles $\{\z_1, \ldots, \z_{N_p}\}$ in the $z$ plane. With the poles
located, the least square problem
\[
\sum_{j=1}^{N_p} \frac{r_j}{\z_j-z} \approx g(z)
\]
computes the residues $\{r_j\}$.

\subsection{Implementation}\label{sec:num}

To implement this algorithm numerically, we need to take care several numerical issues.

\begin{itemize}

\item The semi-infinite matrix in \eqref{eq:lsR} and \eqref{eq:lsL}. In the implementation, we pick
  a value $\dmax$ that is believed to be the upper bound of the number of poles and form the matrix
  \begin{equation}
    H = 
    \begin{bmatrix}
      \hg_{-1} & \hg_{-2} & \cdots & \hg_{-\dmax} \\
      \hg_{-2} & \hg_{-3} & \cdots & \hg_{-(\dmax+1)} \\
      \vdots & \vdots & \vdots & \vdots\\
      \hg_{-l} & \hg_{-(l+1)} & \cdots & \hg_{-(\dmax+l-1)}
    \end{bmatrix}
    \quad\text{or}\quad
    H =
    \begin{bmatrix}
      \hg_{1} & \hg_2 & \cdots & \hg_{\dmax} \\
      \hg_{2} & \hg_3 & \cdots & \hg_{(\dmax+1)} \\
      \vdots & \vdots & \vdots & \vdots\\
      \hg_{l} & \hg_{l+1} & \cdots & \hg_{(\dmax+l-1)}
    \end{bmatrix},
    \label{eq:H}
  \end{equation}
  respectively for \eqref{eq:lsR} and \eqref{eq:lsL}, with $l$ satisfying $l \ge \dmax$. We find that
  in practice $l=\dmax$ is enough.

\item Numerical estimation of the rank $d$ in \eqref{eq:lsR} and \eqref{eq:lsL}. To address this,
  let $s_1,s_2,\ldots,s_{\dmax}$ be the singular values of the matrix $H$. The numerical rank is
  chosen to be the smallest $d$ such that $s_{d+1}/s_1$ is below the noise level.
  
\item Computation of the vector $p$. We first compute the singular value decomposition (SVD) of
\[
\begin{bmatrix}
  \hg_{-1} & \hg_{-2} & \cdots & \hg_{-(d+1)} \\
  \hg_{-2} & \hg_{-3} & \cdots & \hg_{-(d+2)} \\
  \vdots & \vdots & \vdots & \vdots\\
  \hg_{-l} & \hg_{-(l+1)} & \cdots & \hg_{-(d+l)}
\end{bmatrix}
\;\text{or}\;
\begin{bmatrix}
  \hg_{1} & \hg_2 & \cdots & \hg_{d+1} \\
  \hg_{2} & \hg_3 & \cdots & \hg_{d+2} \\
  \vdots & \vdots & \vdots & \vdots\\
  \hg_{l} & \hg_{l+1} & \cdots & \hg_{d+l}
\end{bmatrix},
\]
respectively for \eqref{eq:lsR} and \eqref{eq:lsL}. $p$ is then chosen to be the last column of the
$V$ matrix.

\item The matrix $H$ in \eqref{eq:H} requires the Fourier transform $\hg_k$ from $k=-(\dmax+l-1)$ to
  $(\dmax+l-1)$. In the random access model, we choose an even $N_s \ge 2(\dmax+l)$ and define for
  $n=0,\ldots,N_s-1$
  \begin{equation}
    t_n = \exp\left(i \frac{2\pi n}{N_s}\right), \quad z_n = -\sqrt{ab}\frac{t_n+1}{t_n-1}.
    \label{eq:tzloc}
  \end{equation}
  Using samples $\{g(t_n)\}$ at the points $\{t_n\}$ corresponds to approximating
  \eqref{eq:integral} with the trapezoidal rule. The trapezoidal rule is exponentially convergent
  for smooth functions when the step size $h = \frac{2\pi}{N_s}$ is sufficient small. In the current
  setting, this corresponds to
  \[
  h \ll \sqrt{\frac{a}{b}}, \quad\text{i.e.,}\quad
  N_s \gg \sqrt{\frac{b}{a}}.
  \]
  Applying the fast Fourier transform to $\{g(t_n)\}$ gives $\{\hg_k\}$ for
  $k=-\frac{N_s}{2},\ldots,\frac{N_s}{2}-1$. Among them, $\hg_{-(\dmax+l-1)},\ldots,\hg_{(\dmax+l-1)}$
  are used to form the $H$ matrix in \eqref{eq:H}.
  
\item In the Matsubara model, $g(z)$ is only given at the Matsubara grid
  \[
  z_n =
  \begin{cases}
    2n \frac{\pi i}{\beta},      & \text{for bosons},\\
    (2n+1) \frac{\pi i}{\beta},  & \text{for fermions}.
  \end{cases}
  \]
  computing the integral \eqref{eq:Fourier} is not convenient in the $t$ space since the images $t_n
  = \frac{z_n-\sqrt{ab}}{z_n+\sqrt{ab}}$ are not uniformly distributed. Instead, using
  \eqref{eq:mobius} the integral is equal to
  \[
  \frac{1}{2\pi i}\int_{+i\infty}^{-i\infty} g(z) \left(\frac{z-\sqrt{ab}}{z+\sqrt{ab}}\right)^{-(k+1)}     \frac{2\sqrt{ab}}{(z+\sqrt{ab})^2} \dz
  \approx \frac{-1}{\beta} \sum_{n\in\Z} g(z_n)  \left(\frac{z_n-\sqrt{ab}}{z_n+\sqrt{ab}}\right)^{-(k+1)} \frac{2\sqrt{ab}}{(z_n+\sqrt{ab})^2},
  \]
  in the $z$ variable, where the last step uses the trapezoidal quadrature on the Matsubara
  grid. The trapezoidal rule is exponentially convergent in the regime $a\gg\pi/\beta$. Since the
  last sum is over all integers, it needs to be truncated between $-N_m$ and $N_m$ for some integer
  $N_m$. Noticing that the terms in the sum decays only quadratically, $N_m$ is typically chosen to
  be quite large for a good accuracy.

\item The least square solve for $\{r_j\}$. Using the $z_n$ points in \eqref{eq:tzloc}, we solve the
  following system
  \[
  r = \argmin_{x\in\C^{N_p}} \frac{1}{2}\|A x - b\|^2,\quad
  A = \left[ \frac{1}{\z_j-z_n} \right]_{n,j},\quad
  b = \begin{bmatrix}
    g(z_1) \\
    \ldots\\
    g(z_{N_s})
  \end{bmatrix},
  \]
  The entries of $r$ are the residues $\{r_j\}$.

\end{itemize}

\subsection{Matrix-valued version}

Let us comment on the matrix-valued version \eqref{eq:matver}. The algorithm remains essentially the
same. Below we list the differences.
\begin{itemize}
\item $\hG_k$ is now the matrix-valued Fourier coefficients from the samples $G(t_n)\equiv
  G(z(t_n))$.

\item The SVD is applied to the $l N_b^2 \times (d+1)$ matrices
\[
\begin{bmatrix}
  \cv(\hG_{-1}) & \cv(\hG_{-2}) & \cdots & \cv(\hG_{-(d+1)}) \\
  \cv(\hG_{-2}) & \cv(\hG_{-3}) & \cdots & \cv(\hG_{-(d+2)}) \\
  \vdots & \vdots & \vdots & \vdots\\
  \cv(\hG_{-l}) & \cv(\hG_{-(l+1)}) & \cdots & \cv(\hG_{-(d+l)})
\end{bmatrix}
\;\text{or}\;
\begin{bmatrix}
  \cv(\hG_{1}) & \cv(\hG_2) & \cdots & \cv(\hG_{d+1}) \\
  \cv(\hG_{2}) & \cv(\hG_3) & \cdots & \cv(\hG_{d+2}) \\
  \vdots & \vdots & \vdots & \vdots\\
  \cv(\hG_{l}) & \cv(\hG_{l+1}) & \cdots & \cv(\hG_{d+l})
\end{bmatrix}
\]
where $\cv(\cdot)$ turns a matrix into a column vector.

\item The least square problem is applied to
\[
R = \argmin_{X\in\C^{N_p\times N_b^2}} \frac{1}{2}\|A X - B\|^2, \quad
A = \left[ \frac{1}{\z_j-z_n} \right]_{n,j},\quad
B = \begin{bmatrix}
  \rv(G(z_1)) \\
  \ldots\\
  \rv(G(z_{N_s}))
\end{bmatrix},
\]
where $\rv(\cdot)$ turns a matrix into a row vector. Each row of $R$ is then reshaped back to the
$N_b\times N_b$ matrix $R_j$. In the special case of $R_j = v_j v_j^*$, $v_j$ can be further
constructed by applying a rank-1 approximation to $R_j$.

\end{itemize}

\subsection{Special cases and extensions}\label{sec:ext}

Below we include a few comments concerning special cases and direct extensions.
\begin{itemize}
\item We have assumed that the poles reside in the two disks in the $z$ plane. In many applications,
  it is known that the poles are actually on the real axis. In such as case, the Fourier
  coefficients $\hg_k$ and hence the matrix $H$ are real. Therefore, a real SVD can be used while
  determining the rank $d$ and the coefficients $(p_0,\ldots,p_d)$. Finally, the roots of $p(z)$ are
  also real. These considerations can significantly improve stability as shown in Section
  \ref{sec:res}.


  
\item We have not specified any noise model. If the noise model is known, it is possible to denoise
  the values $g(z_n)$ before applying the algorithm described. Such a denoising step can potentially
  improve the accuracy and stability of pole locations.

\item The algorithm can also be extended to the general setting, where the imaginary axis $i\R$ is
  replaced with any simple curve in the Riemann sphere. If the curve is smooth, the extension is
  straightforward as the trapezoidal quadrature can still be applied. When the curve is non-smooth,
  special quadrature is often needed for good accuracy.
\end{itemize}

\section{Numerical results}\label{sec:res}


This section presents a few numerical examples. In all examples, $a=1$, $b=100$. The
noise added to $g(z)$ is multiplicative:
\[
g_{\text{noisy}} = g_{\text{exact}} \cdot ( 1 + \sigma N_{\C}(0,1) ).
\]
This is a reasonable model since in many applications the magnitude of the noise is often
proportional to the magnitude of the signal. For each example, we present the numerical results for
both the random access model and the Matsubara model. For the random access model, $N_s=1024$. 
For the Matsubara model, $N_m=10^6$ and $\beta=10\pi$.

\begin{figure}[h!]
  \begin{tabular}{cc}
    \includegraphics[scale=0.3]{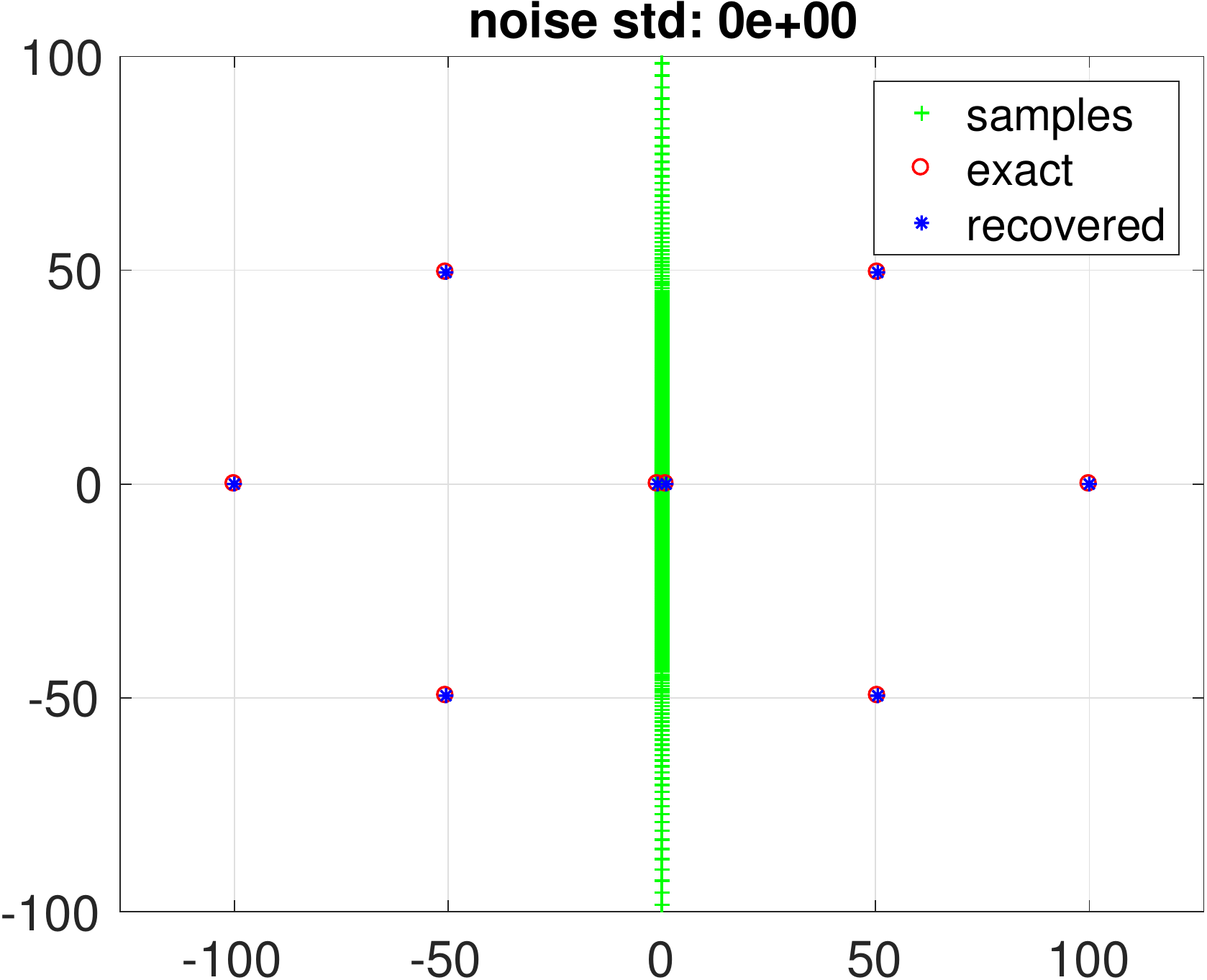} & \includegraphics[scale=0.3]{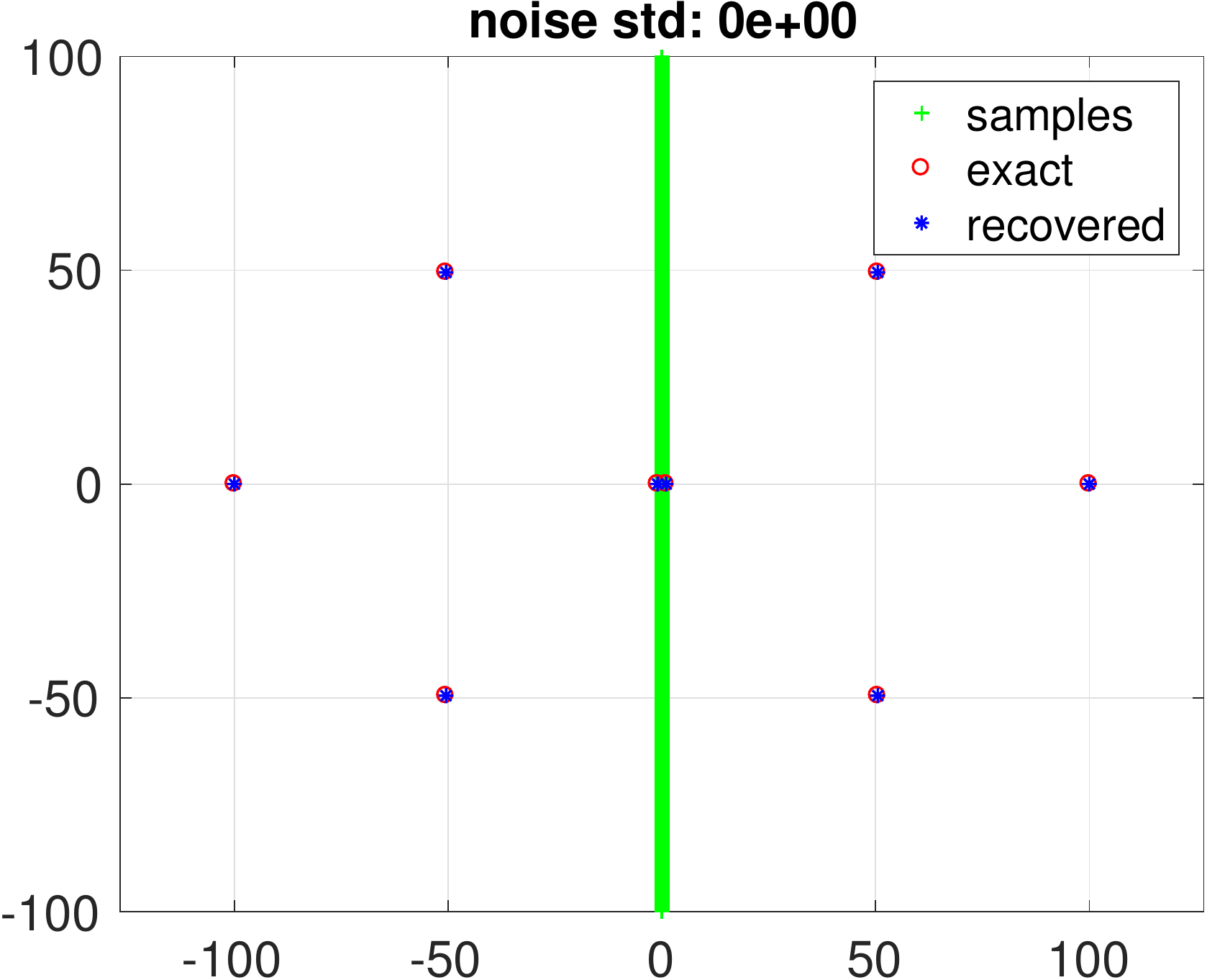}\\
    \includegraphics[scale=0.3]{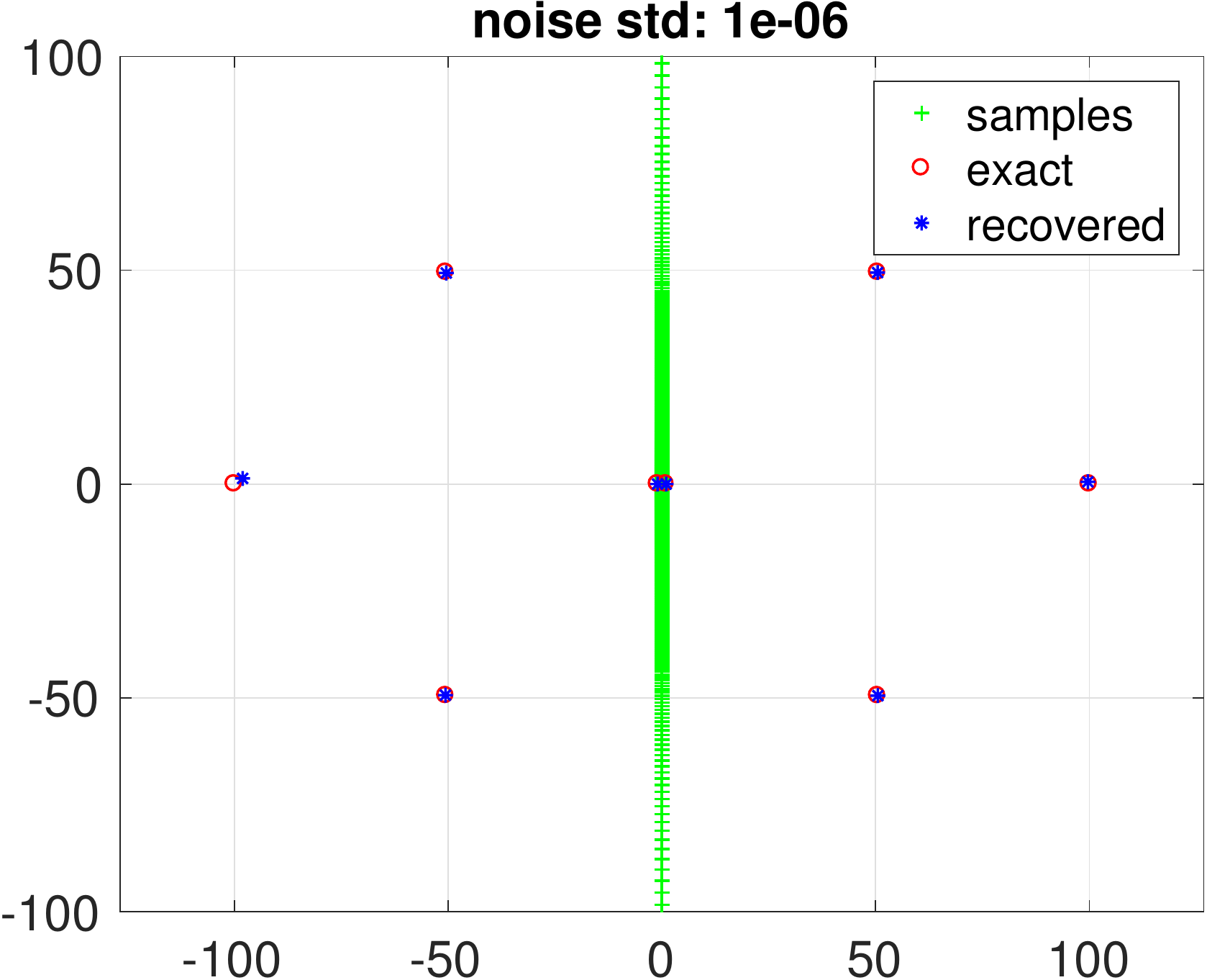} & \includegraphics[scale=0.3]{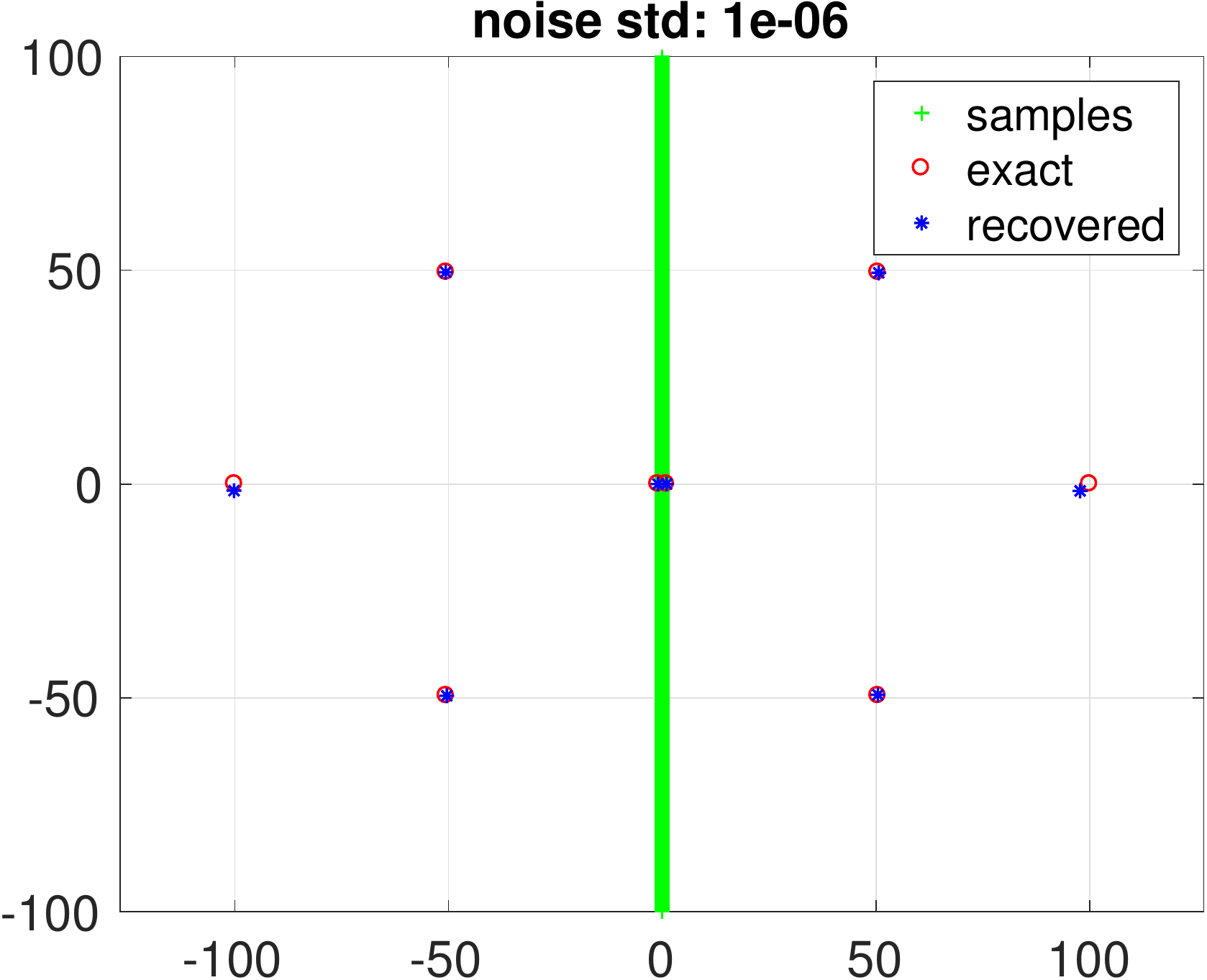}\\
    \includegraphics[scale=0.3]{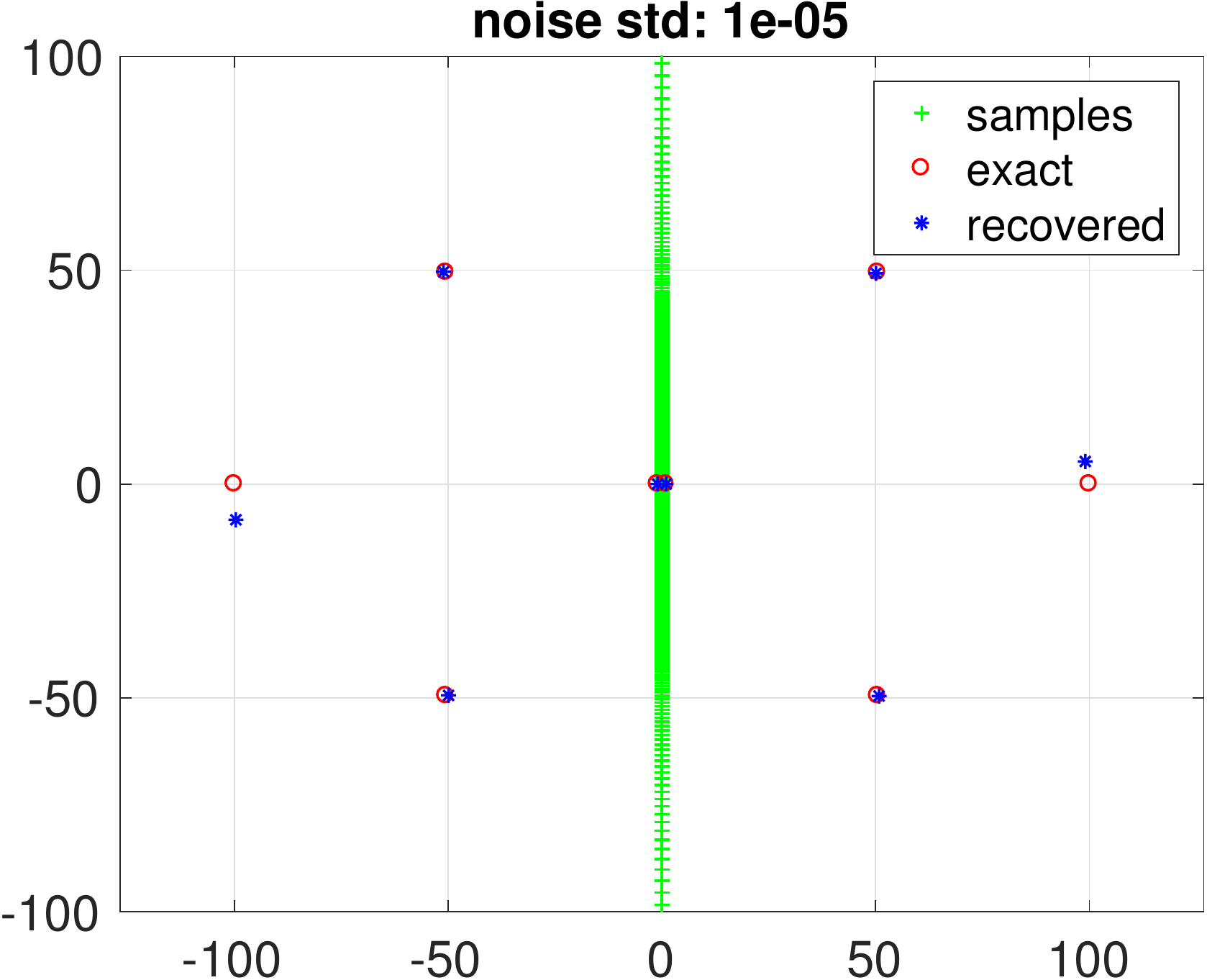} & \includegraphics[scale=0.3]{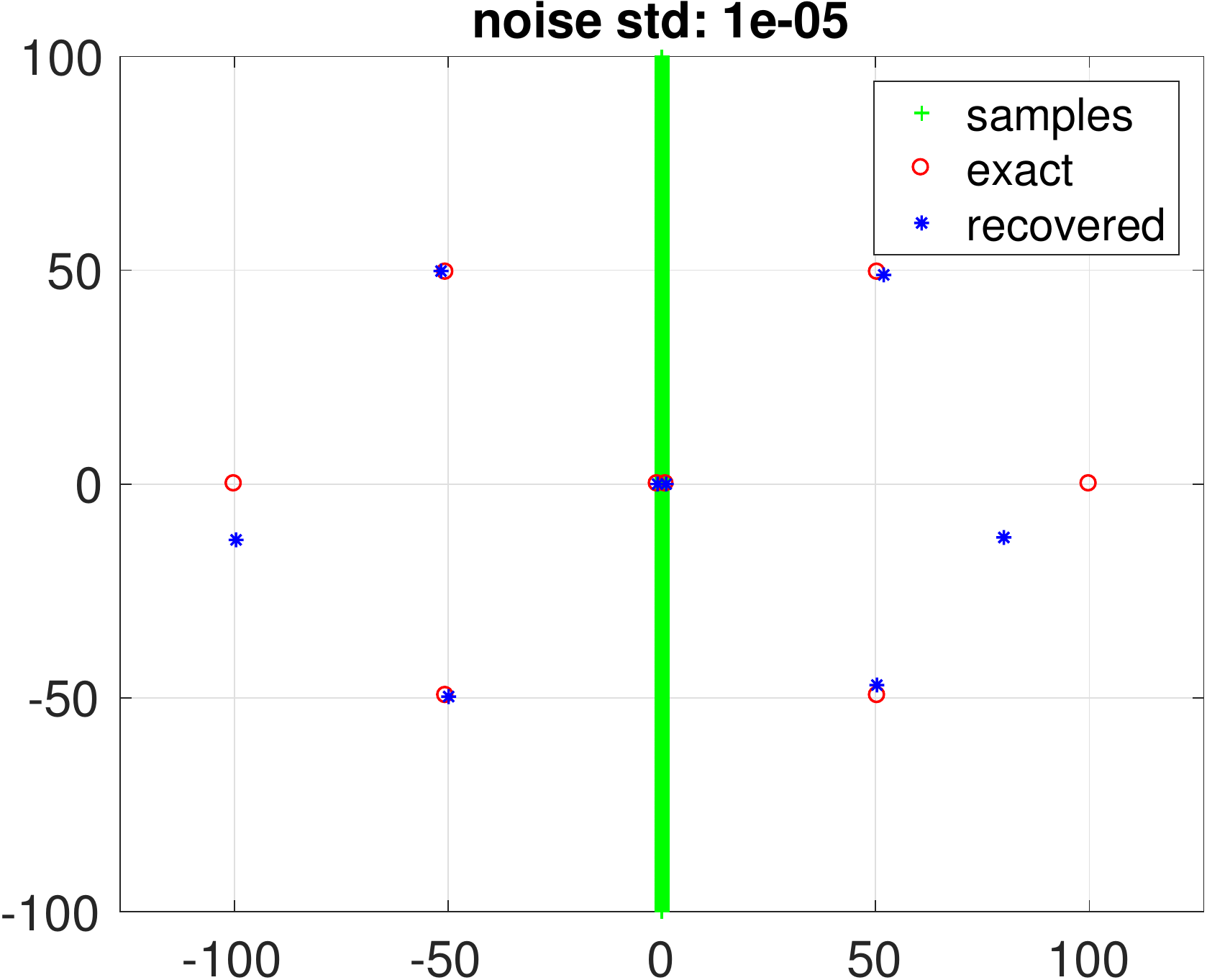}\\
    \includegraphics[scale=0.3]{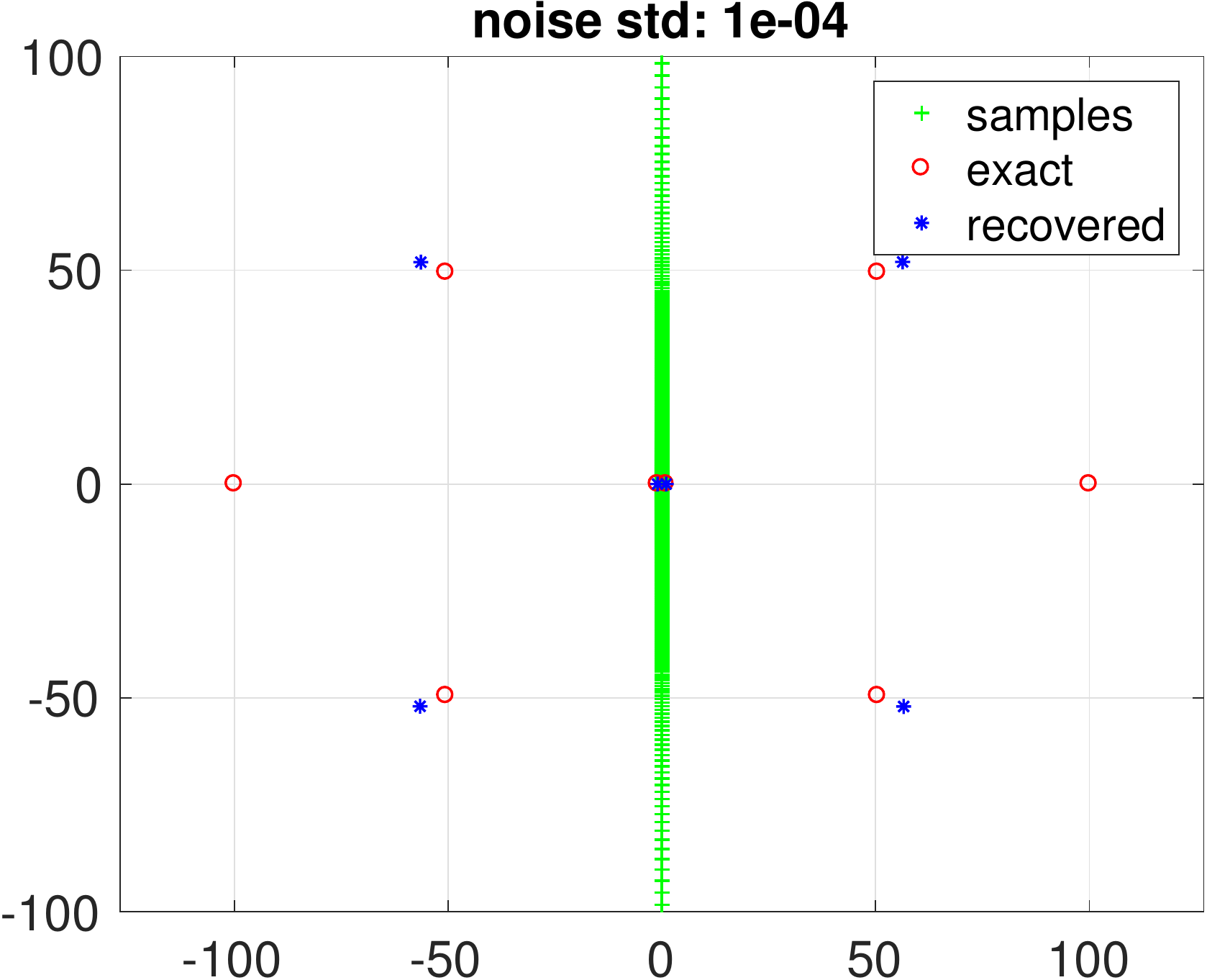} & \includegraphics[scale=0.3]{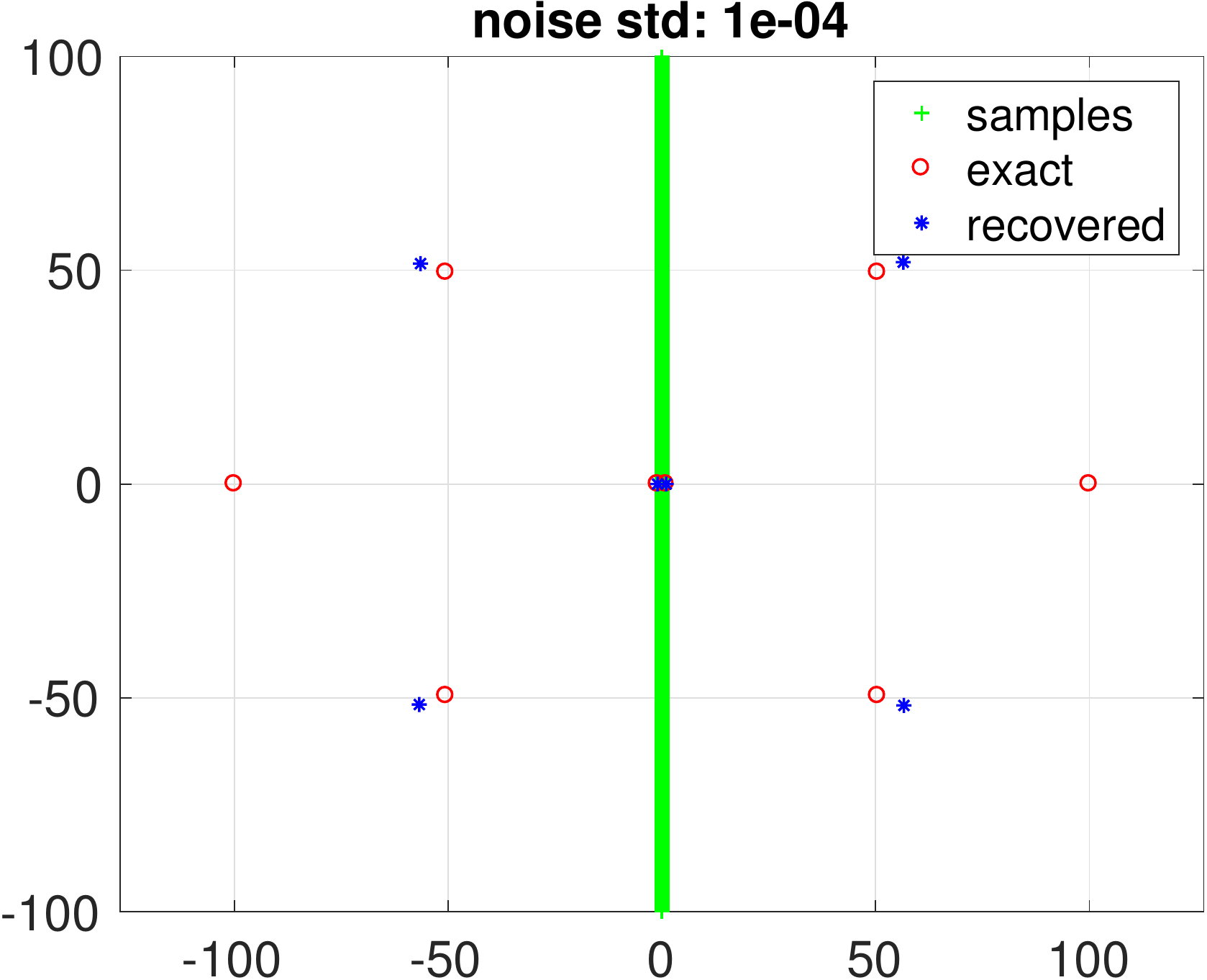}\\
    Random access model & Matsubara model
  \end{tabular}
  \caption{Complex pole locations, with different levels of noise. Left: random access model. Right:
    Matsubara model.}
  \label{fig:cplx}
\end{figure}

{\bf Example 1.} We first consider the case of complex pole locations. Within each circle, we place
four poles and the residues $\{r_j\}$ are of unit order. Figure \ref{fig:cplx} plots the results at
the noise level $\sigma=0$, $10^{-6}$, $10^{-5}$, and $10^{-4}$, where the left and right columns
are for the random access and Matsubara models, respectively. The results show that
\begin{itemize}
\item At $\sigma=0$, the algorithm gives perfect reconstruction at machine accuracy.
\item At $\sigma=10^{-6}$, the poles are accurately identified.
\item At $\sigma=10^{-5}$, the number of poles are correctly recovered, though the locations of the
  two poles far from $i\R$ are wrong.
\item At $\sigma=10^{-4}$, only the six poles close to $i\R$ are identified.
\end{itemize}

\begin{figure}[h!]
  \begin{tabular}{cc}
    \includegraphics[scale=0.3]{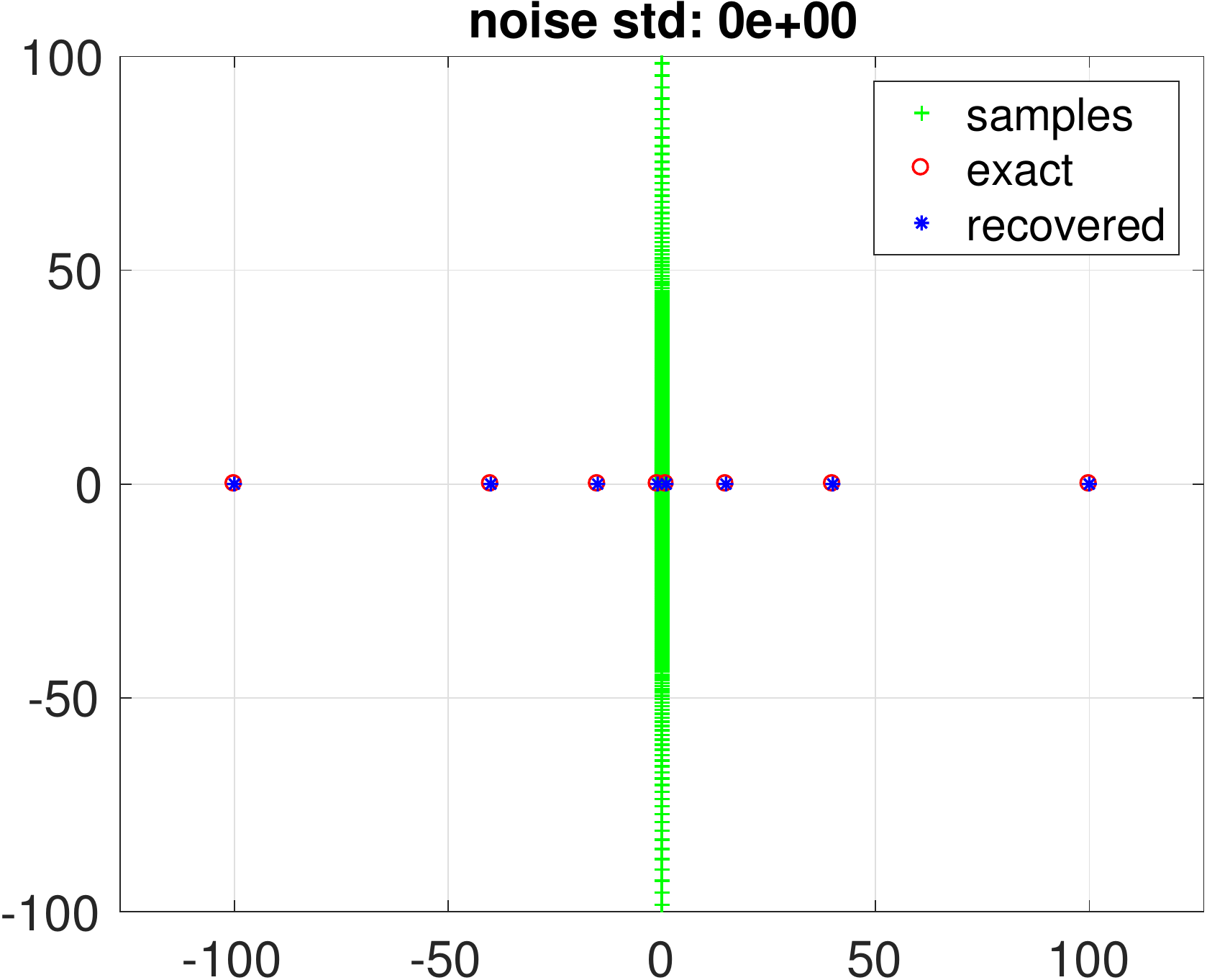} & \includegraphics[scale=0.3]{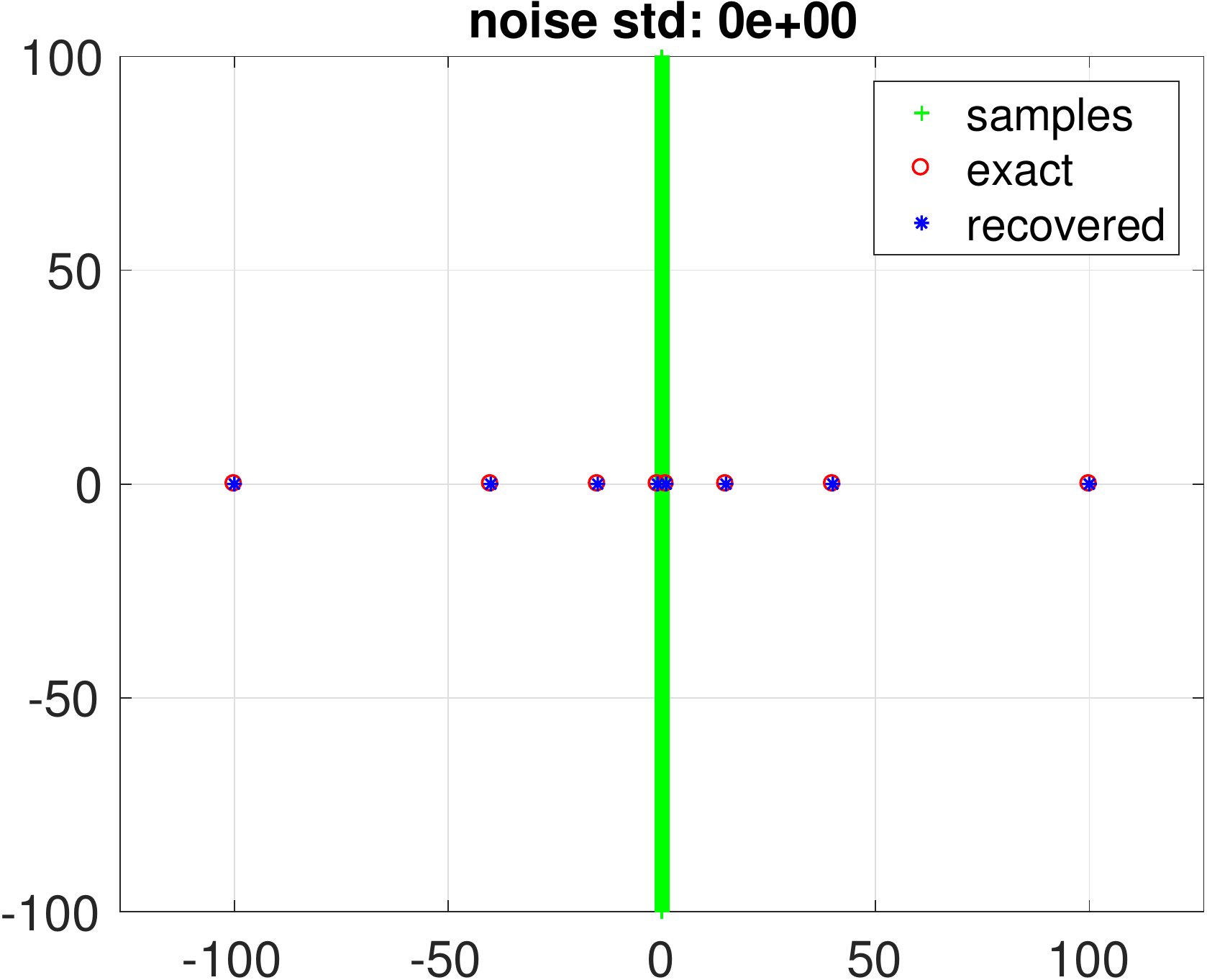}\\
    \includegraphics[scale=0.3]{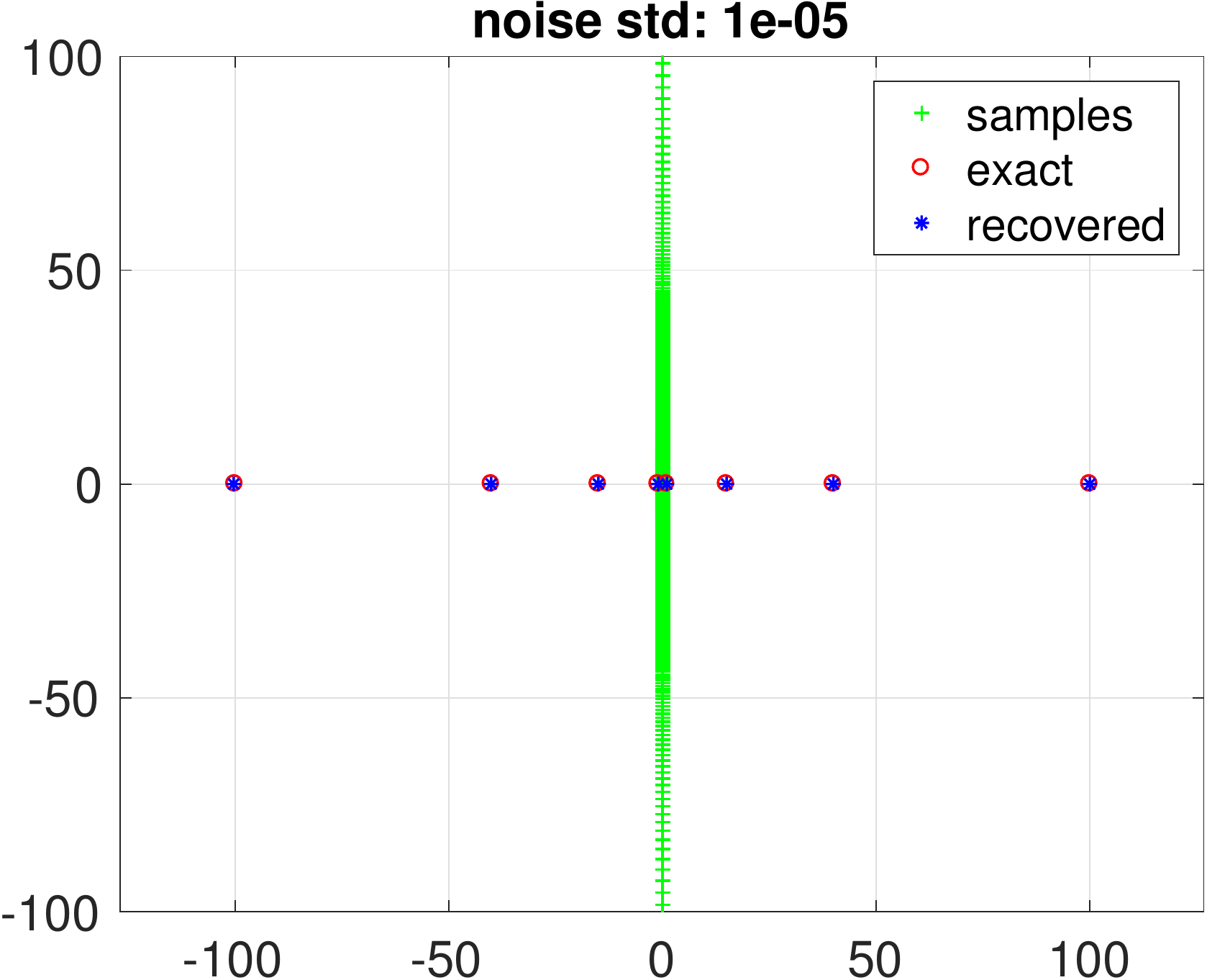} & \includegraphics[scale=0.3]{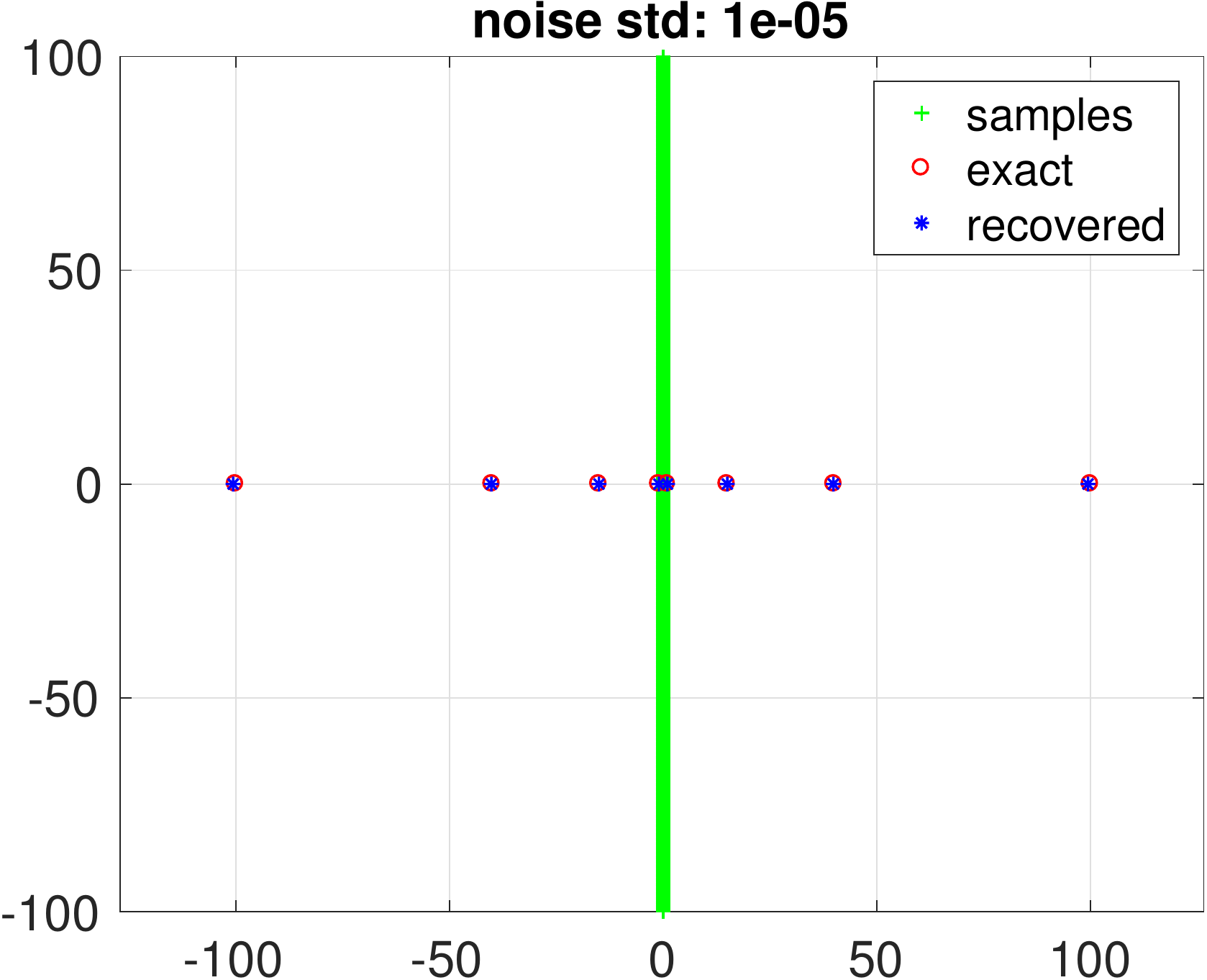}\\
    \includegraphics[scale=0.3]{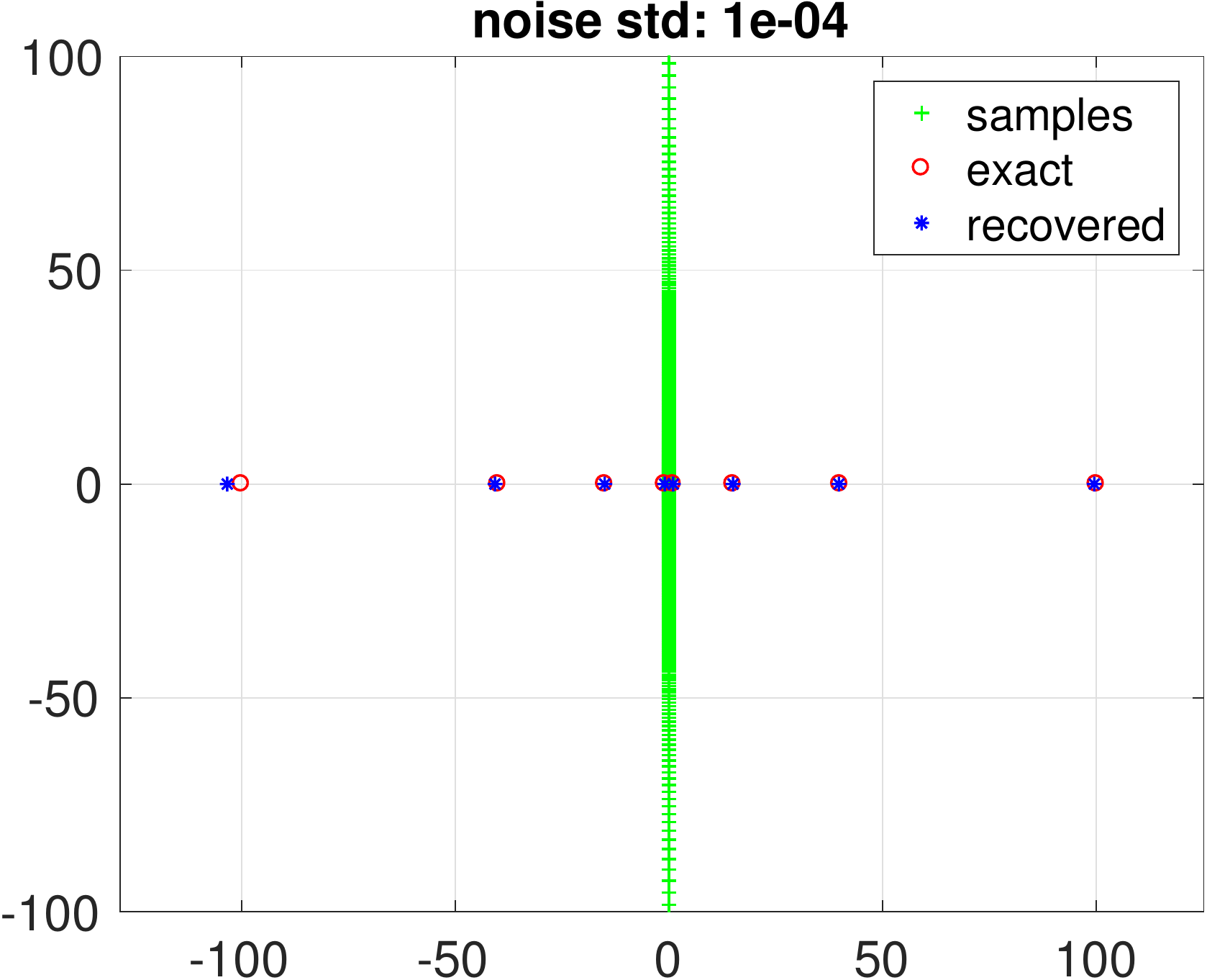} & \includegraphics[scale=0.3]{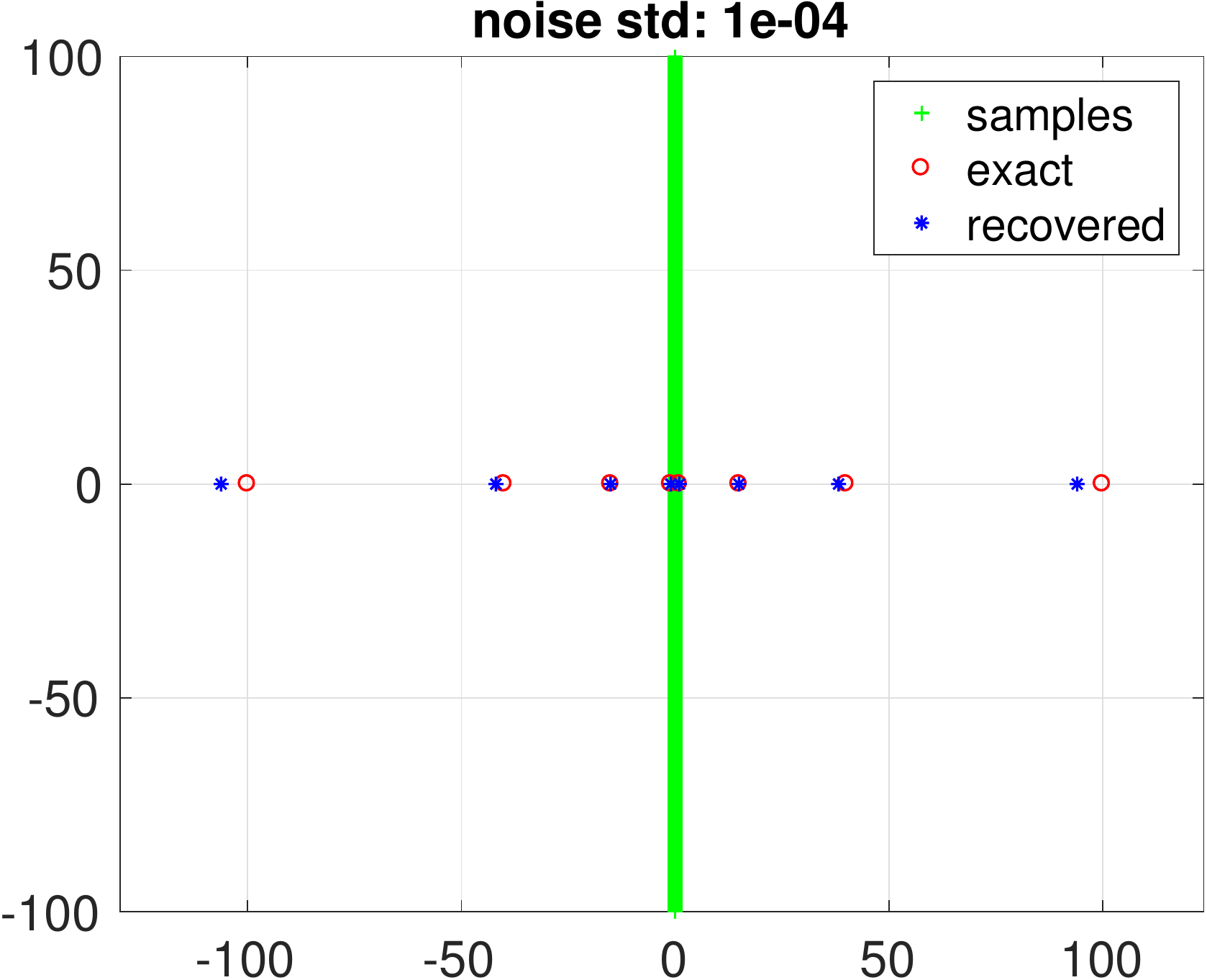}\\
    \includegraphics[scale=0.3]{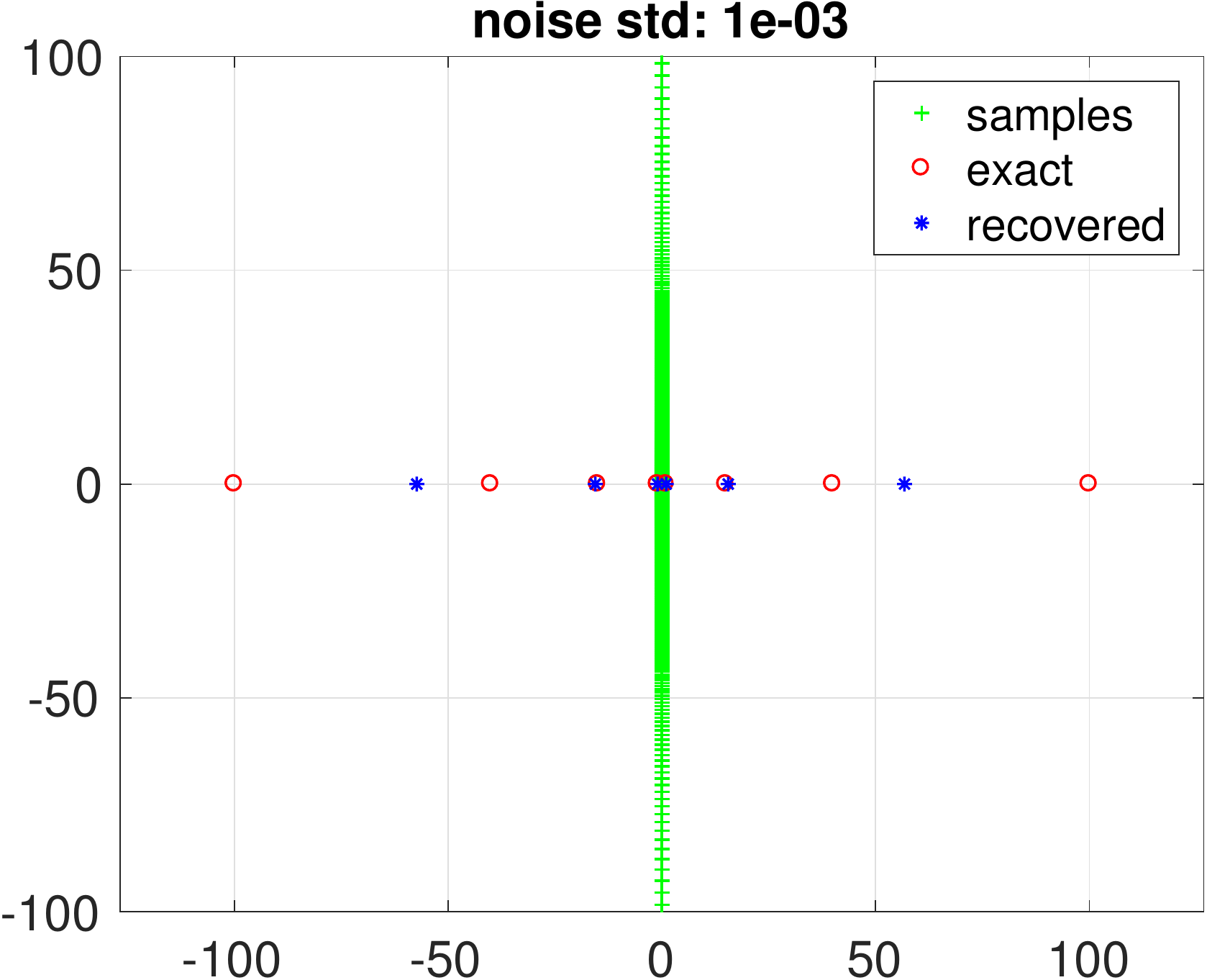} & \includegraphics[scale=0.3]{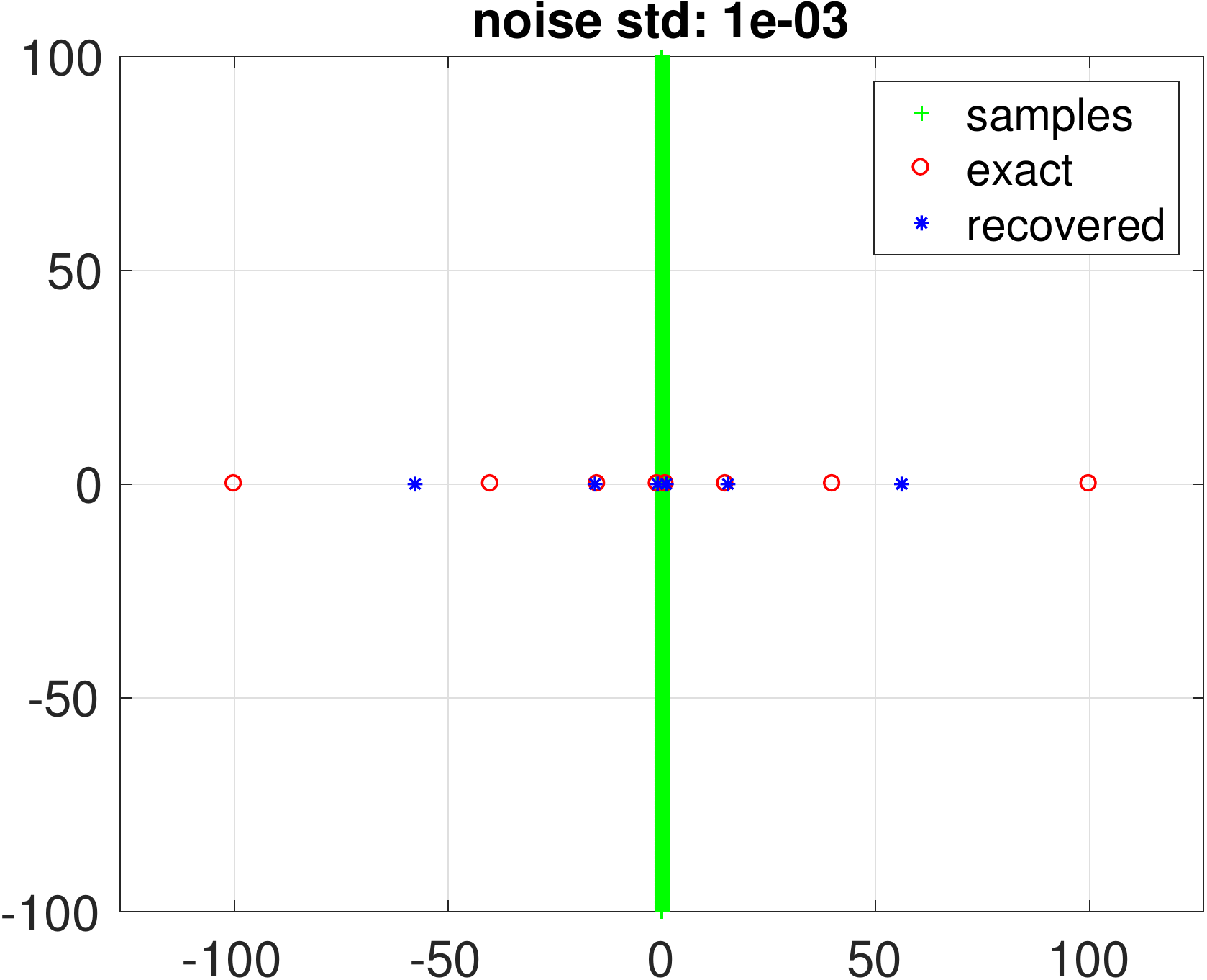}\\
    Random access model & Matsubara model
  \end{tabular}
  \caption{Real pole locations, with different levels of noise. Left: random access model. Right:
    Matsubara model.}
  \label{fig:real}
\end{figure}

{\bf Example 2}. Next we consider the case of real pole locations. Within each circle, there are 4
poles and the residues $\{r_j\}$ are again of unit order. Figure \ref{fig:real} summarizes the
results at the noise level $\sigma=0$, $10^{-5}$, $10^{-4}$, and $10^{-3}$.
\begin{itemize}
\item At $\sigma=0$, the algorithm gives perfect reconstruction.
\item At $\sigma=10^{-5}$, the poles are also recovered perfectly.
\item At $\sigma=10^{-4}$, the pole locations are recovered accurately, though with some errors for
  the two poles farthest away from $i\R$.
\item At $\sigma=10^{-3}$, only the six poles close to $i\R$ are identified.
\end{itemize}
A comparison with the previous example suggests that enforcing the real constraints significantly
improves the stability and accuracy of the algorithm.

\begin{figure}[h!]
  \begin{tabular}{cc}
    \includegraphics[scale=0.3]{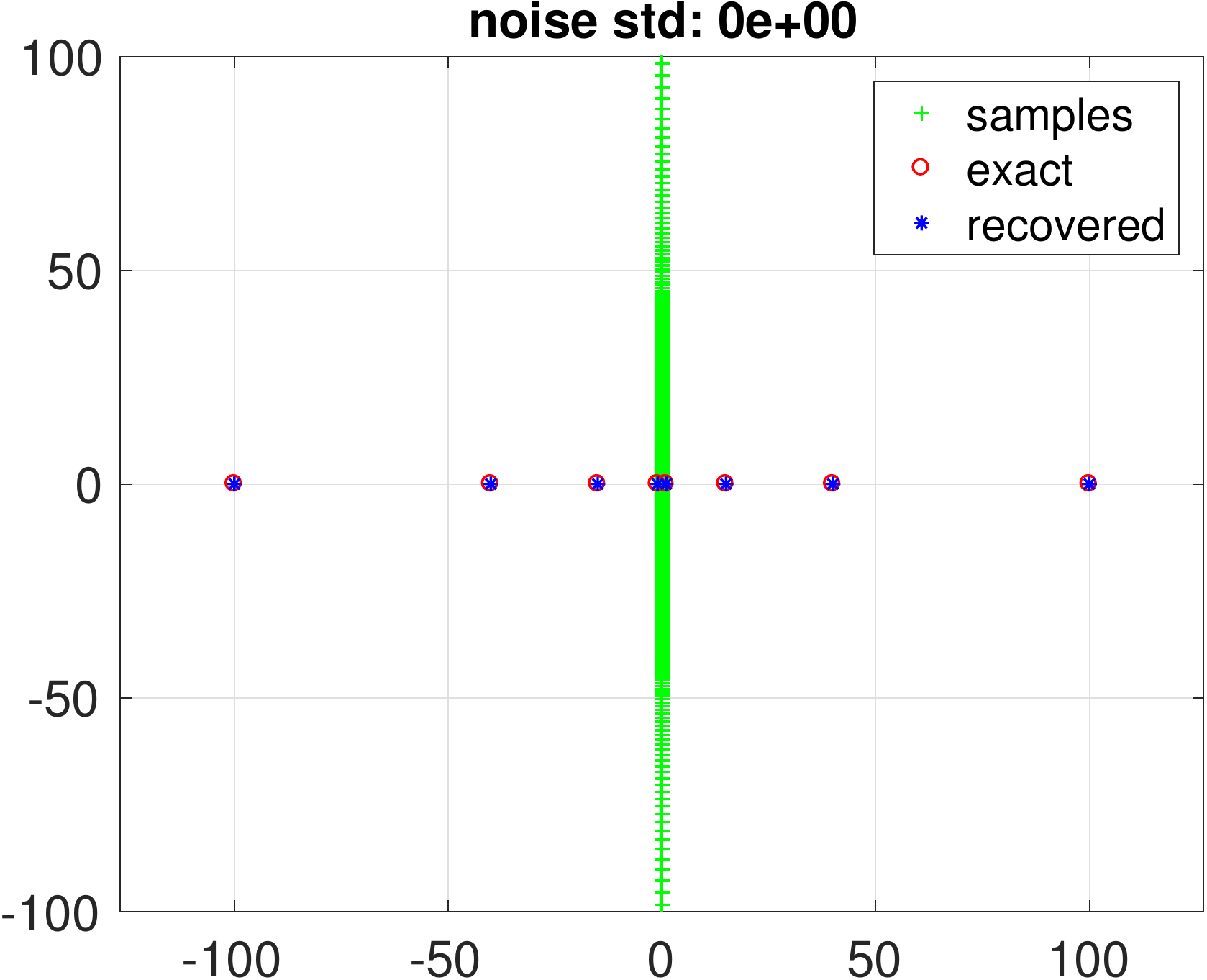} & \includegraphics[scale=0.3]{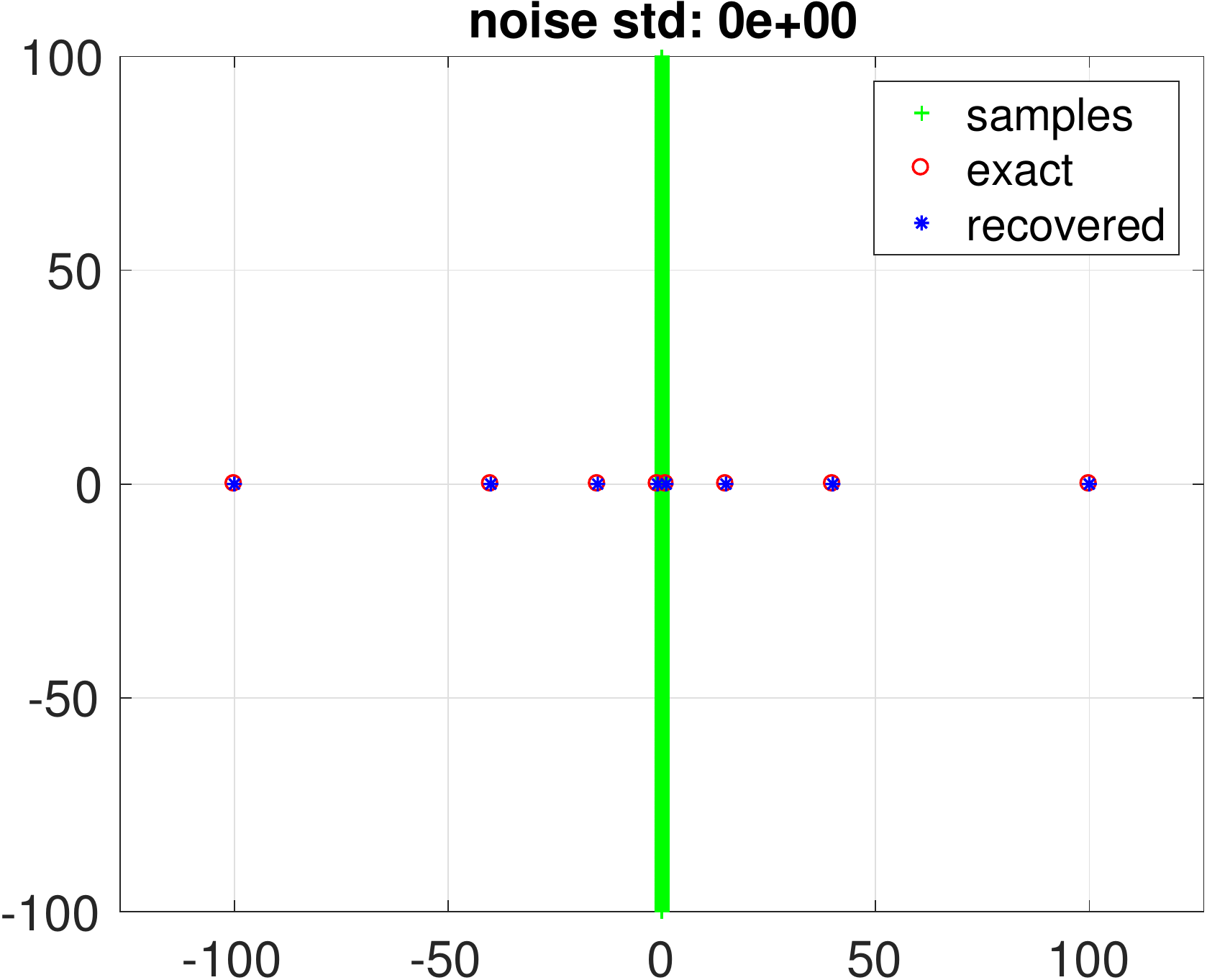}\\
    \includegraphics[scale=0.3]{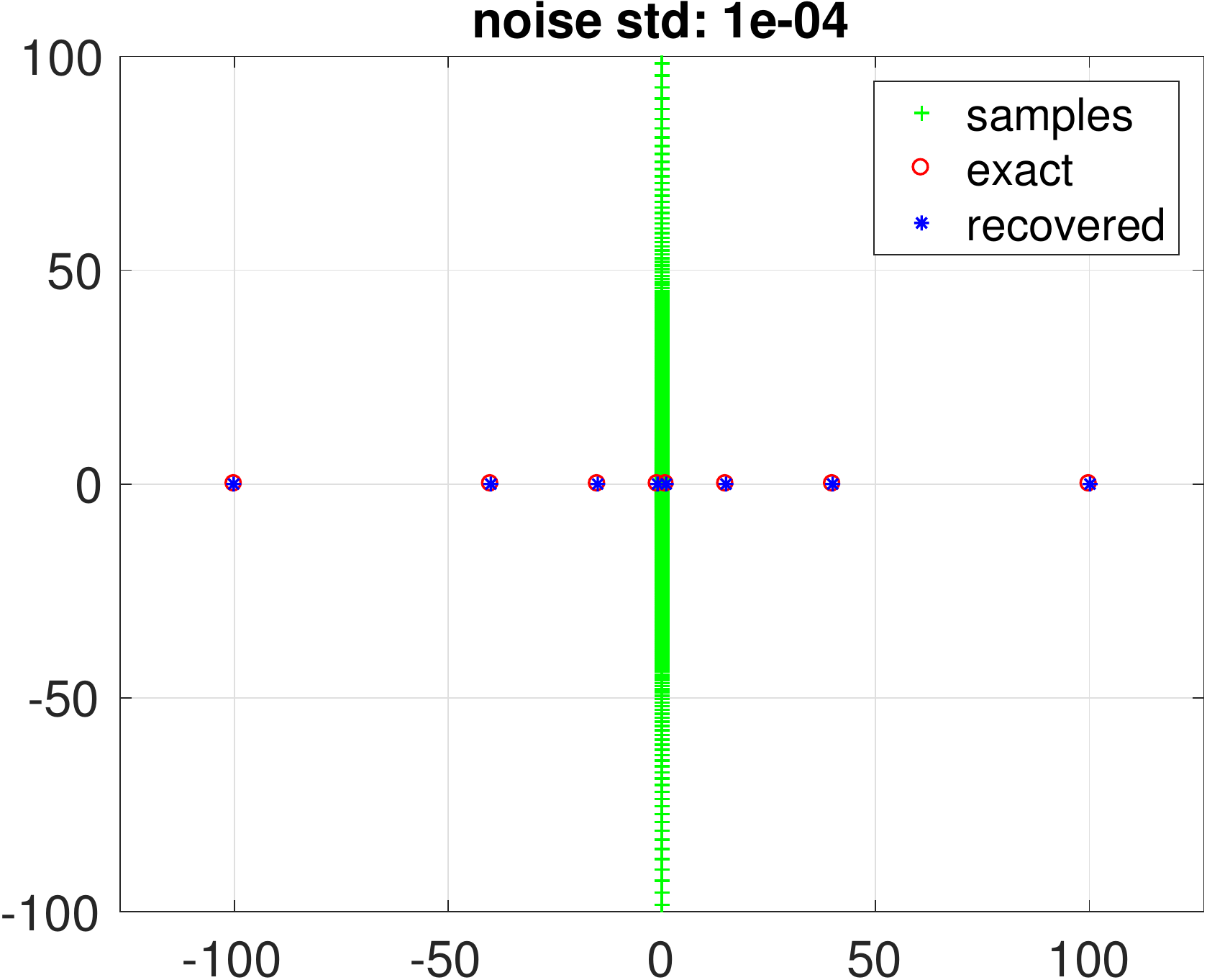} & \includegraphics[scale=0.3]{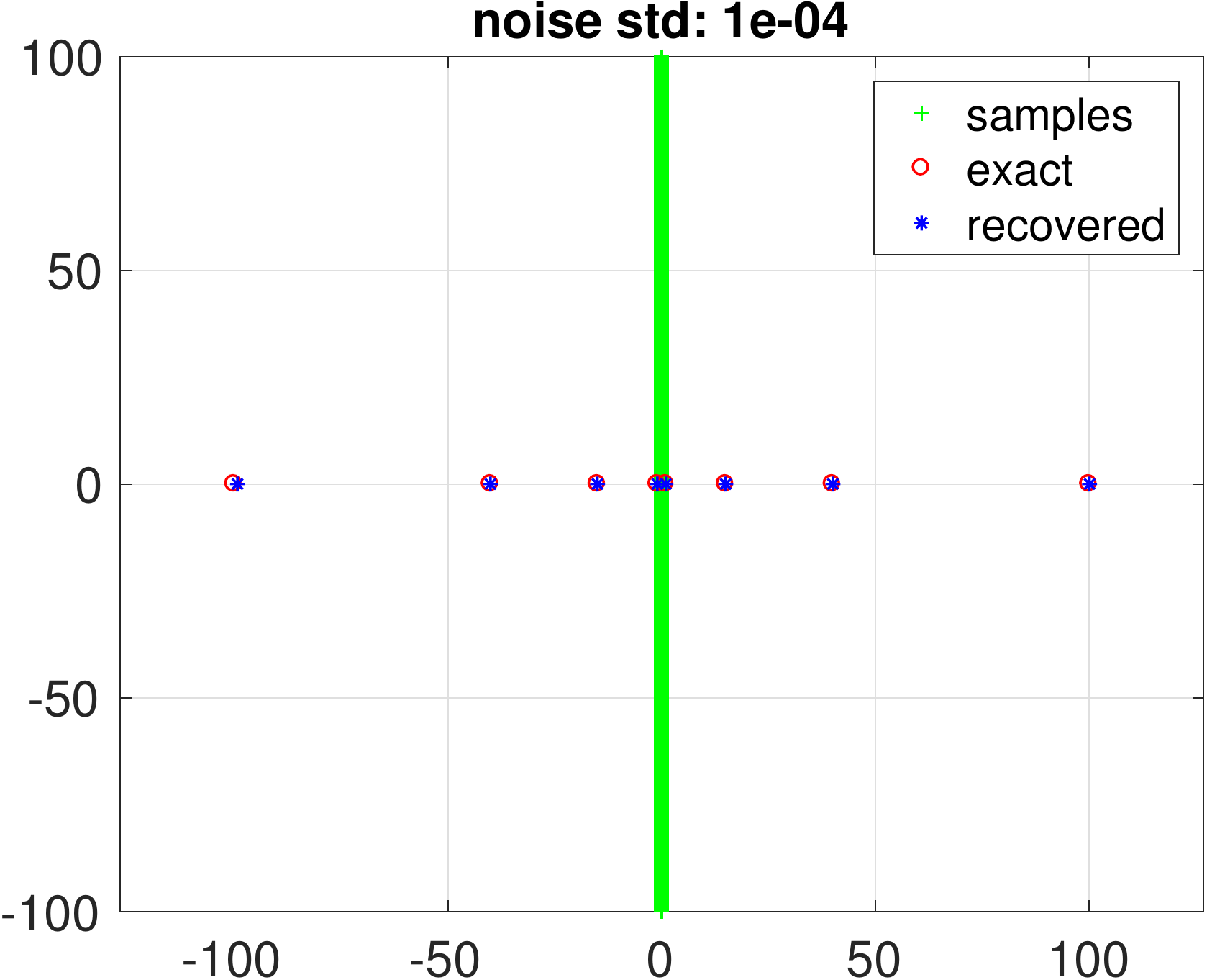}\\
    \includegraphics[scale=0.3]{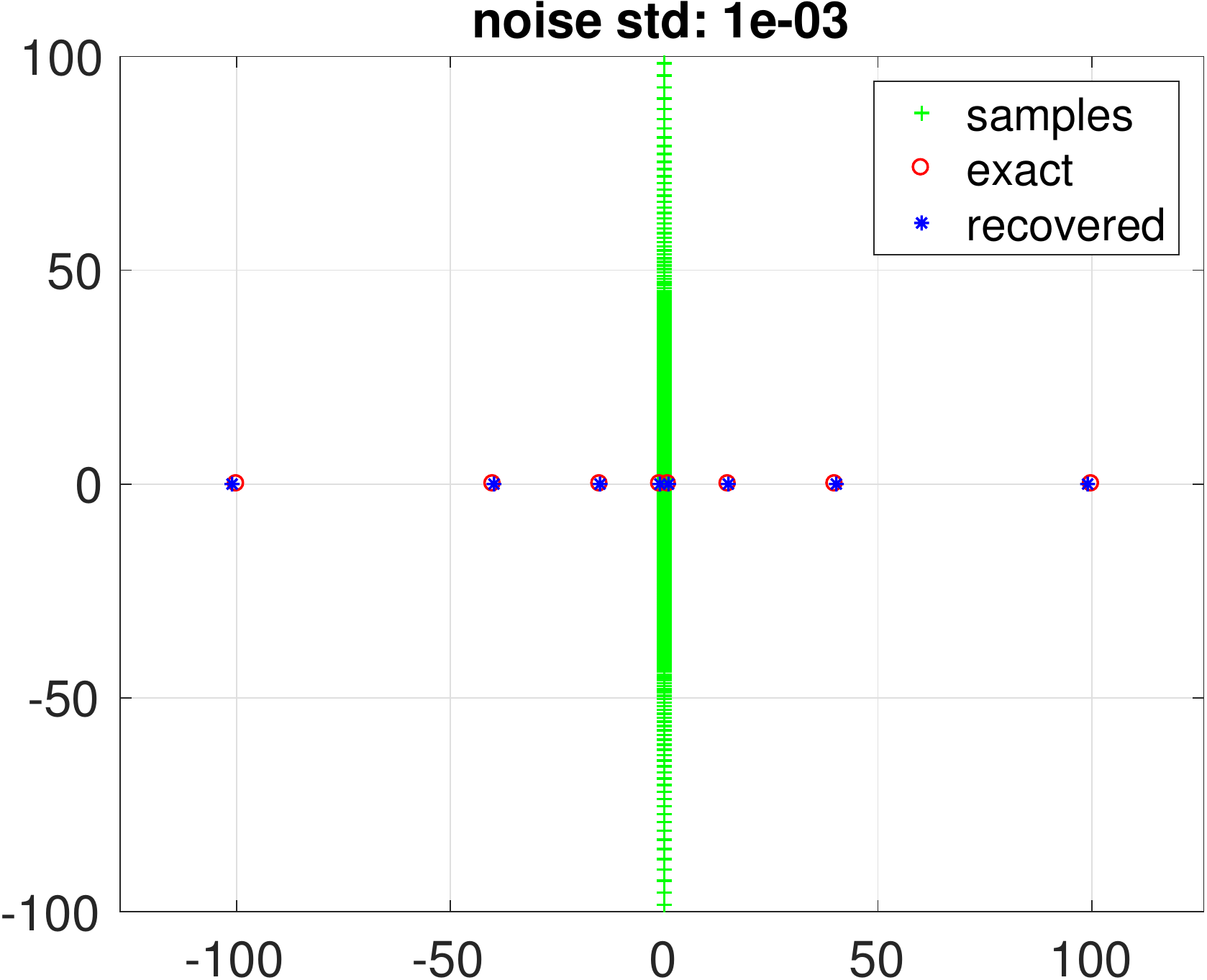} & \includegraphics[scale=0.3]{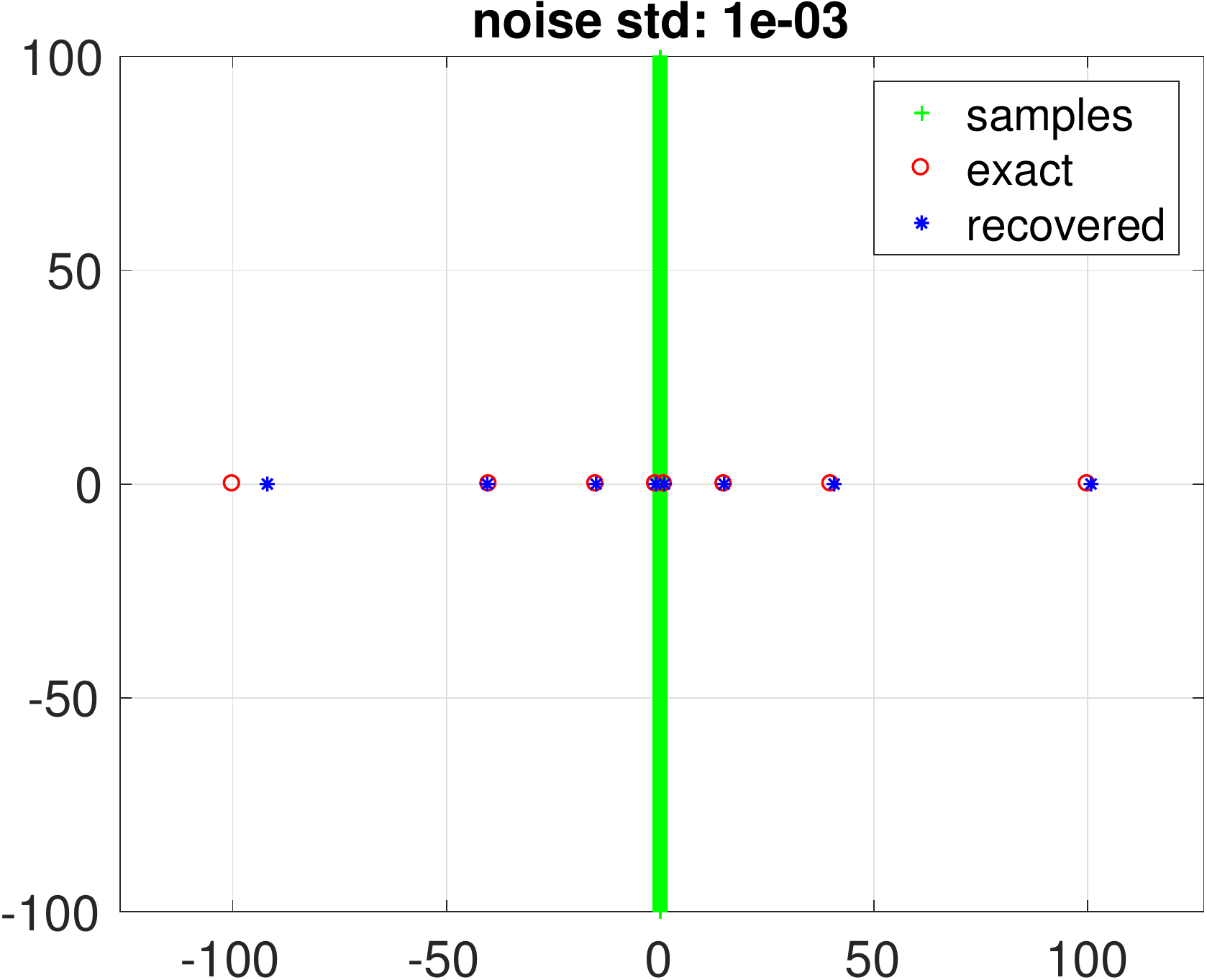}\\
    \includegraphics[scale=0.3]{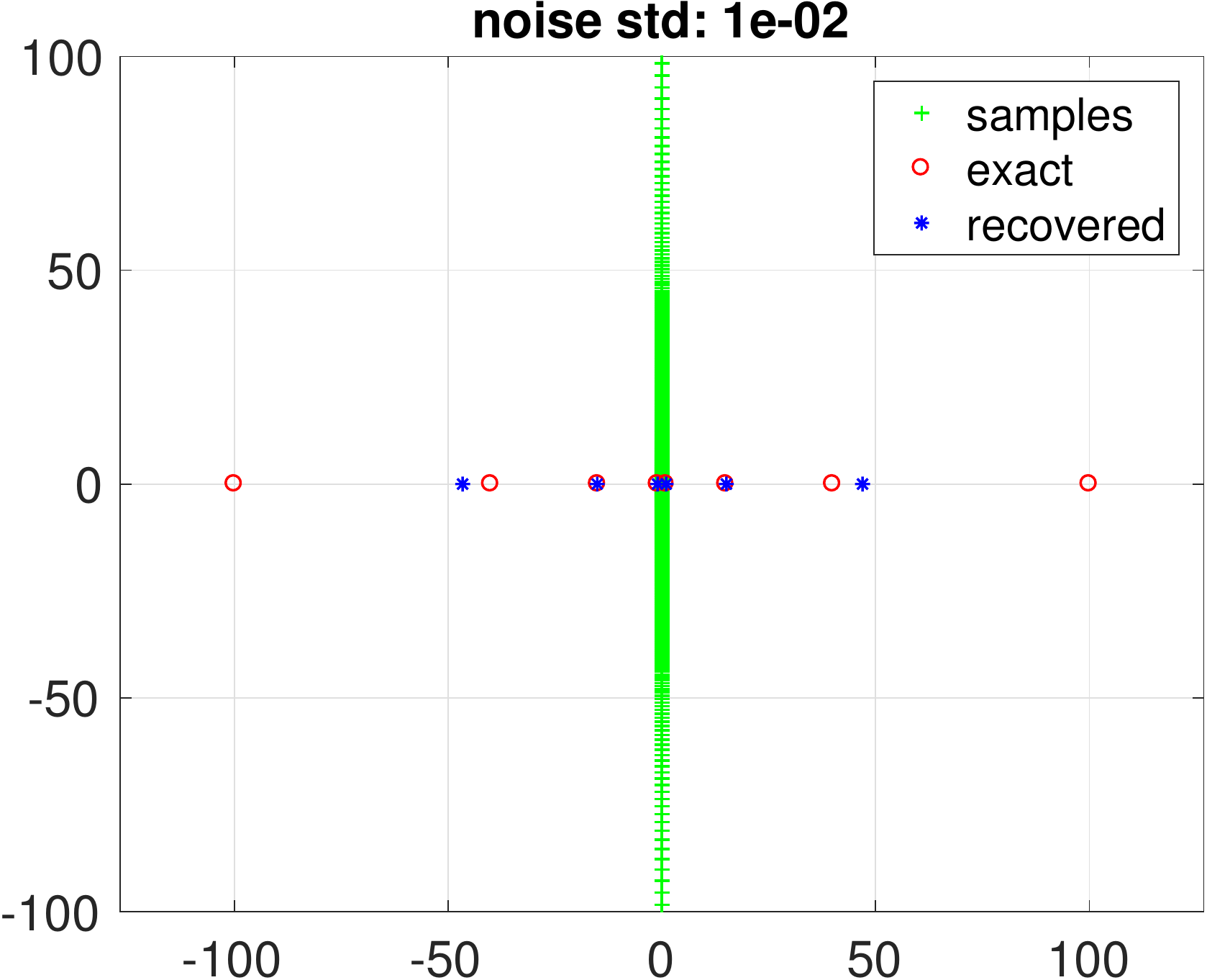} & \includegraphics[scale=0.3]{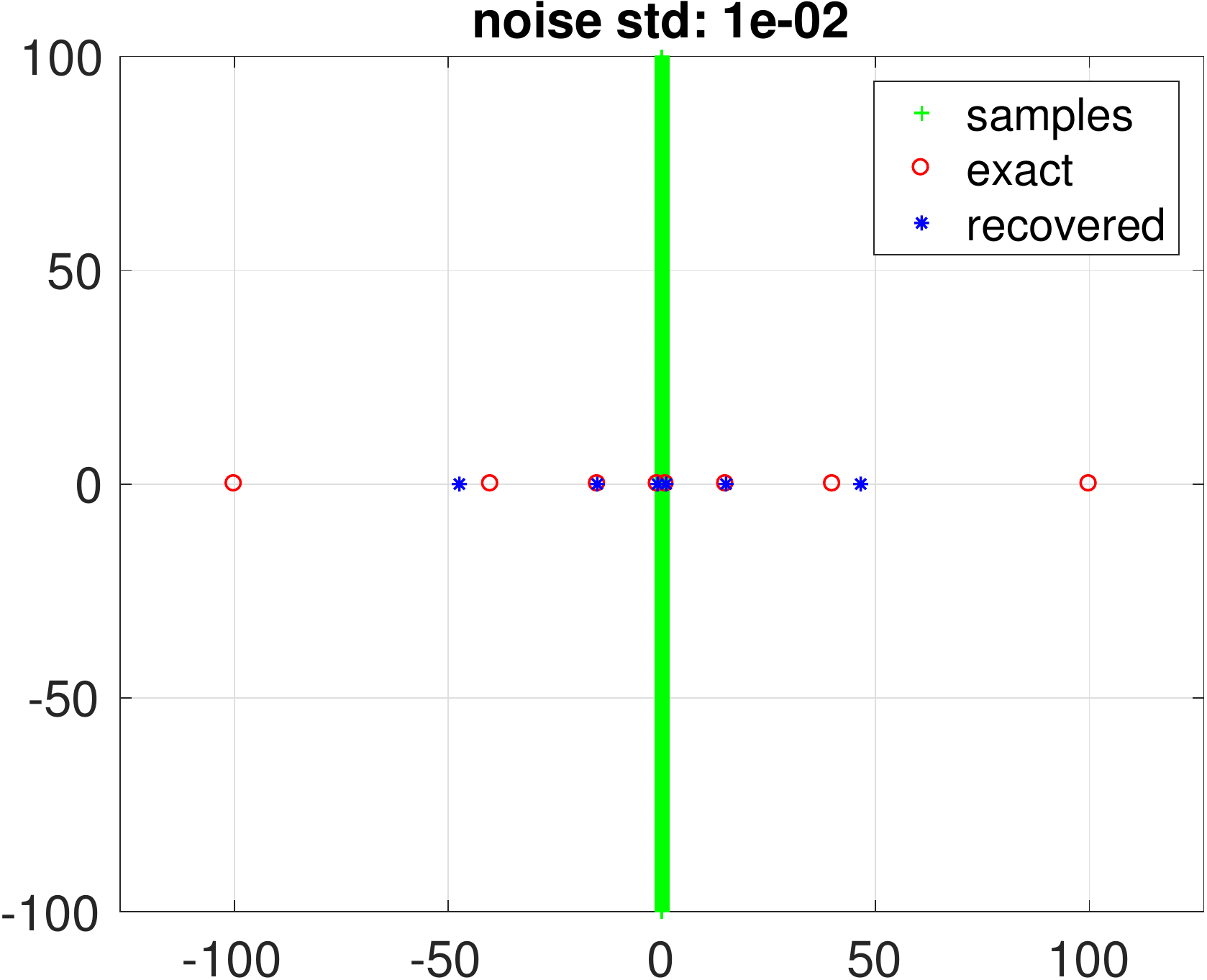}\\
    Random access model & Matsubara model
  \end{tabular}
  \caption{Matrix case with real poles, with different levels of noise. Left: random access
    model. Right: Matsubara model.}
  \label{fig:mtrx}
\end{figure}

{\bf Example 3}. Finally, we consider the matrix-valued version. The dimension $N_b$ of the matrix
$R_j$ is set to be $N_b=4$. When other parameters are fixed, larger values of $N_b$ significantly
improve the accuracy since it effectively provides more data. Within each circle, there are again 4
poles and the residues $\{v_j\}$ (and equivalently $\{R_j\}$) are of unit order. Figure
\ref{fig:mtrx} summarizes the results at the noise level $\sigma=0$, $10^{-4}$, $10^{-3}$,
and $10^{-2}$.

\begin{itemize}
\item At $\sigma=0$, the algorithm again gives perfect reconstruction.
\item At $\sigma=10^{-4}$, the reconstruction is near perfect.
\item At $\sigma=10^{-3}$, the pole locations are recovered with good accuracy, though there are
  some errors for the two poles away from $i\R$.
\item At $\sigma=10^{-2}$, only the six poles close to $i\R$ are identified.
\end{itemize}
Noticing that the noise level in this example is much higher than the ones used in the previous
examples, the results confirm that the matrix-valued version is easier, especially when $N_b$ is
large.

\bibliographystyle{abbrv}

\bibliography{ref}

\end{document}